\documentclass{article}

\begin{document}

\title{K\"ahler-Einstein metrics on Fano manifolds, II: limits with cone angle less than $2\pi$.}
\author{Xiuxiong Chen, Simon Donaldson and Song Sun\footnote{X-X. Chen was partly supported by National Science Foundation grant; S. Donaldson and S. Sun are partly supported by the European Research Council award No 247331.}}
\date{\today}
\maketitle


\newtheorem{thm}{Theorem}
\newtheorem{prop}{Proposition}
\newtheorem{lem}{Lemma}
\newtheorem{cor}{Corollary}
\newtheorem{defn}{Definition}
\newcommand{\tomega}{\tilde{\omega}}
\newcommand{\tD}{\tilde{D}}
\newcommand{\Ric}{{\rm Ricci}}
\newcommand{\bC}{{\bf C}}
\newcommand{\bQ}{{\bf Q}}
\newcommand{\bP}{{\bf P}}
\newcommand{\bZ}{{\bf Z}}
\newcommand{\bR}{{\bf R}}
\newcommand{\db}{\overline{\partial}}
\newcommand{\cD}{{\cal D}}
\newcommand{\uf}{\underline{f}}
\newcommand{\us}{\underline{s}}
\newcommand{\uv}{\underline{v}}
\newcommand{\dbd}{i \partial \overline{\partial}}
\newcommand{\Euc}{{\rm Euc}}

\section{Introduction}

\

This is the second of a series of three papers which provide proofs of results announced in \cite{kn:CDS0}. We review the background briefly but refer the reader also to the introduction to \cite{kn:CDS1}. 
Let $X$ be a Fano manifold of complex dimension $n$. Let $\lambda\geq 0$ be an integer and $D$ be a smooth divisor in the linear system $\vert -\lambda K_{X}\vert$. (For purely expository purposes, we allow the case when $\lambda=0$,
in which case $D$ is the empty set.)  For $\beta\in (0,1]$ we consider a K\"ahler-Einstein metric $\omega$ with  a cone angle $2\pi\beta$ along $D$.
To be precise, we take this to mean that when $\beta<1$ the metric is defined by a potential in the class ${\cal C}^{2,\alpha,\beta}$ defined in \cite{kn:Dona11}, for all $\alpha<\beta^{-1}-1$. We recall this definition in a little more detail in Section 3.2 below. Once we have a metric of this type, higher regularity results  have been given by Jeffres, Mazzeo and Rubinstein in \cite{kn:JMR}; but for our purposes we do not need to use those results. When $\beta=1$ we mean that the metric is a smooth K\"ahler-Einstein metric in the usual sense, but in the present paper we always consider $\beta<1$.
  The Ricci curvature of such a K\"ahler-Einstein  metric $\omega$, with cone singularities, is $(1-\lambda(1-\beta))\omega$. For our purposes here we can  suppose that  $\beta> 1-\lambda^{-1}$ so that $\omega$ has strictly positive Ricci curvature.  The point of this paper is to study the convergence properties of sequences of such metrics.
 Thus we consider a sequence of triples $(X_{i}, D_{i},\omega_{i})$ with fixed $\lambda$ and fixed Hilbert polynomial $\chi(X, K_{X}^{-p})$ and a sequence of cone angles $\beta_{i}$ with limit $\beta_{\infty}$. (For our main application we can take $(X_{i}, D_{i})=(X, D)$, so only the cone angle varies, but even in that case it becomes notationally easier to write $(X_{i}, D_{i})$.)
We will also usually denote the positive line bundle  $K_{X}^{-1}$ by $L$, partly to streamline notation and partly because many of our results would apply to more general polarised manifolds. In this paper we consider the case when the limit $\beta_{\infty}$ is strictly less than 1. In the sequel we will consider the case when $\beta_{\infty}=1$ and also explain, in more detail than in \cite{kn:CDS0}, how our results lead to the main theorem announced there.

To state our main result we need to recall some theory of K\"ahler-Einstein metrics on singular varieties. A general reference for this is \cite{kn:EGZ}. If $W$ is a normal variety we write $K_{W}$ for the canonical line bundle over the smooth part $W_{0}$ of $W$.
\begin{defn} A $\bQ$-Fano variety is a normal projective variety $W$ which \begin{itemize} \item is $\bQ$-Gorenstein,  so some power $K_{W}^{m_{0}}$ extends to a line bundle over $W$;
\item for some multiple $m$ of $m_{0}$, is embedded in projective space by sections of $K_{W}^{-m}$;
\item  has Kawamata log terminal (KLT) singularities. \end{itemize}
\end{defn}

(For theoretical purposes we do not need to distinguish between $m$ and $m_{0}$, but in reality one expects that $m_{0}$ can often be taken substantially smaller than $m$.)

\begin{defn}
Let $W\subset \bC\bP^{N}$ be a $\bQ$-Fano variety as above, embedded by sections of $K_{W}^{-m}$. A {\it K\"ahler metric on $W$ in the class $c_{1}$} is a K\"ahler metric $\omega$ on the smooth part $W_{0}$ of $W$ which has the form
$ m^{-1} \omega_{FS}+ \dbd \phi$ where $\omega_{FS}$ is the restriction of the Fubini-Study metric and $\phi$ is a continuous function on $W$, smooth on $W_{0}$.
\end{defn}

Another way of expressing this uses metrics on line bundles. We get a reference metric on the line bundle $K_{W}^{-m}$ from the identification with ${\cal O}(1)\vert_{ W}$. This gives a reference metric
$h_{0}$ on $K_{W}^{-1}$ over $W_{0}$ and a K\"ahler form $\omega$ is defined by a metric $h=e^{-\phi} h_{0}$ on $K_{W}^{-1}$. Thus we will often write $\omega=\omega_{h}$. Of course we could change $h$ by a non-zero multiple  without changing $\omega_{h}$. We say that a metric $h$ on the line bundle $K_{W}^{-1}$ over $W_{0}$ is continuous on $W$ if the induced metric on $K_{W}^{-m}$ extends continuously to $W$.

A metric $h$ on $K_{W}^{-1}$ defines a volume form $\Omega_{h}$ on $W_{0}$: if $\sigma$ is an element of the fibre over a point with $\vert \sigma\vert=1$ then the volume form at that point is $(\sigma\wedge\overline{\sigma})^{-1}$. 
\begin{defn}
Let $W$ be a $\bQ$-Fano variety and $\omega_{h}$ be a K\"ahler metric on $W$ in the class $c_{1}$. The metric is K\"ahler-Einstein if 
\begin{equation} \omega_{h}^{n}=  \Omega_{h} \end{equation}
on the smooth set $W_{0}\subset W $.
\end{defn}

Now fix $\lambda\geq 0$ as above.  Let $s$ be a section of the line bundle $K_{W}^{-\lambda}$ over $W_{0}$ such that
\begin{itemize}
\item if $\lambda=0$, the section $s$ is a non-zero constant;
\item if $\lambda>0$, the power $s^{m}$ extends to  non-trivial  holomorphic section of $K_{W}^{-m \lambda}$ over $W$.
\end{itemize}
Let $\Delta_{0}$ be the divisor in $W_{0}$ defined by $s$ (so $\Delta_{0}=0$ if $\lambda=0$), and let $\Delta$ be the closure of $\Delta_{0}$ in $W$. Thus $\Delta$ is a Weil divisor in $W$. For $\beta\in (0,1)$ we have an
$\bR$-divisor $(1-\beta)\Delta$.
\begin{defn}
The pair $(W,(1-\beta)\Delta)$ is KLT if for each point $w\in W$ there is a neighbourhood $U$ of $w$, a positive integer $m$, and a non-vanishing holomorphic section $\sigma$ of $K_{W}^{-m}$ over $U\cap W_0$ such that $$\int_{U\cap (W_{0}\setminus {\rm supp} \Delta_{0})} (\sigma \wedge \overline{\sigma})^{-1/m} \vert s\vert_{h}^{2\beta-2} <\infty,$$
where $h$ is any continuous metric on $K_{W}^{-1}$.
\end{defn}
\begin{defn}
Let $W$ be a $\bQ$-Fano variety  and let $\Delta$ be a Weil divisor defined by a section $s$ as above such that $(W, (1-\beta)\Delta)$ is KLT for some $\beta\in (0,1)$.  A {\it weak conical  K\"ahler-Einstein metric} for the triple $(W,\Delta, \beta)$ is a metric $h$ on
$K_{W}^{-1}$ which extends continuously to $W$, which is smooth on $W_{0}\setminus {\rm supp} \Delta_{0}$ and which satisfies the equation
\begin{equation}   \omega_{h}^{n}= \Omega_{h} \vert s \vert_{h}^{2(\beta-1)} \end{equation}
on $W_{0}\setminus {\rm supp} \Delta_{0}$. 
\end{defn}

It would be more precise to say that (2) is an equation for the metric $h$ determined by $\beta$ and the section $s$; however $s$ is determined by $\Delta$ up to a non-zero constant multiple. Given $\Delta$, we can normalise $s$ and the metric $h$, taking non-zero multiples, by requiring that the $L^{2}$ norm of $s$  (computed using $h$ and the volume form $ (n!)^{-1}\omega_{h}^{n}$) is
$1$.

\

With these definitions fixed we can state our main result.
\begin{thm}
Let $X_{i}$ be a sequence of $n$-dimensional Fano manifolds with fixed Hilbert polynomial. Let $D_{i}\subset X_{i}$ be smooth divisors in $\vert -\lambda K_{X_{i}}\vert $, for fixed $\lambda\geq 0$. If $\lambda>0$ let $\beta_{i}\in (0,1)$ be a sequence  converging to a limit $\beta_{\infty}$ with $1-\lambda^{-1}< \beta_{\infty}<1$.  Suppose that there are K\"ahler-Einstein metrics $\omega_{i}$ on $X_{i}$ with cone angle $2\pi\beta_{i}$ along $D_{i}$. Then there is a $\bQ$-Fano variety $W$ and a Weil divisor $\Delta$ in $W$ such that
\begin{enumerate}
\item $(W,(1-\beta_{\infty})\Delta)$ is KLT;
\item  there is a weak conical K\"ahler-Einstein metric $\omega$ for the triple
$(W,\Delta, \beta_{\infty})$;
\item  possibly after passing to a subsequence, there are embeddings $T_{i}:X_{i}\rightarrow \bC\bP^{N}, T_{\infty}:W\rightarrow \bC\bP^{N}$, defined by the complete linear systems $\vert -m K_{X_{i}}\vert$ and $\vert -m K_{W}\vert$ respectively (for a suitable $m$), such that $T_{i}(X_{i})$ converge to $T_{\infty}(W)$ as projective varieties and $T_{i}(D_{i})$ converge to $\Delta$ as algebraic cycles.
\end{enumerate}
\end{thm}
The statement of Theorem 1 is designed to give exactly what we need for the main argument outlined in \cite{kn:CDS0}; no more and no less. One would like to have more information   about the singularities of the limiting metric at points of $\Delta$. This seems quite complicated when the divisor has components of multiplicity greater than $1$ and we leave a discussion for a later paper. When there are no such components then, at smooth points of $\Delta$, we do have complete results. In particular we have

\begin{thm}  With the same hypotheses as in Theorem 1, suppose also that $W$ is smooth and $\Delta\subset W$ is a smooth divisor. Then the metric $\omega$ is a K\"ahler-Einstein metric with cone angle $2\pi \beta_{\infty}$ along $\Delta$.
\end{thm}

\

We now outline the strategy of proof of Theorem 1. By the results of \cite{kn:CDS1} we can approximate each metric $\omega_{i}$ arbitrarily closely in Gromov-Hausdorff distance by smooth K\"ahler metrics with Ricci curvature bounded below by fixed strictly positive number. Since the volume is fixed (as the leading term in the Hilbert polynomial) the metrics satisfy a fixed volume non-collapsing bound.  This means that we can transfer all of the Cheeger-Colding theory of non-collapsed limits of metrics of positive Ricci curvature to our situation.
So, perhaps taking a subsequence, we can suppose that the $X_{i}$ have a Gromov-Hausdorff limit $Z$. (We usually just write $X_{i}$ for $(X_{i}, \omega_{i})$ or, more precisely, the metric space defined by $\omega_{i}$.) The essential problem is to give $Z$ a complex algebraic structure, in the sense of constructing a normal variety $W$ and a homeomorphism $h:Z\rightarrow W$. Then we have to show that $W$ is $\bQ$-Fano,  and that the Gromov-Hausdorff convergence can be mirrored by algebro-geometric convergence in the sense stated in Theorem 1. We also have to show that there is a limiting metric on $W$ which satisfies the appropriate weak K\"ahler-Einstein equation.  In the case when $\lambda=0$, so the divisors do not enter the picture at all, this is what was done in \cite{kn:DS} except for the last statement about the limiting metric (which is not  difficult). So  Theorem 1 is essentially established in \cite{kn:DS} for the case when $\lambda=0$. Of course the real case of interest for us now is when $\lambda>0$, and what we have to do is to adapt the arguments in \cite{kn:DS} to take account of the divisors. The main work takes place in Section 2, establishing that all tangent cones to $Z$ are \lq\lq good''. In Section 3.1 we put together the proof of Theorem 1. Theorem 2 is proved in subsection (3.2).

\section{Structure of Gromov-Hausdorff limits}

\subsection{Preliminaries}
We collect some facts which apply to  rather general length spaces and which will be important to us.  We expect that these are entirely standard results for experts.

\

Let $P$ be a $p$-dimensional length space which is either the based Gromov-Hausdorff limit of a non-collapsed sequence of manifolds with Ricci curvature bounded below, or the cross-section of a tangent cone of a limit of such a sequence.  Let $A$ be a compact subset of $P$. Define a subset ${\cal M}_{A}$ of $(0,\infty)$ as follows. A number $C$ is in ${\cal M}_{A}$ if for all $\epsilon>0$ there is a cover of $A$ by at most $C \epsilon^{2-p}$ balls centred at points of $A$. Define the {\it (codimension 2) Minkowski measure} $m(A)$ to the the infimum of ${\cal M}_{A}$, with the understanding that $m(A)=\infty$ if ${\cal M}_{A}$ is empty. Next we say that $A$ has {\it capacity zero} if for all $\eta>0$ there is a Lipschitz function  $g$, equal to $1$ on a neighbourhood of $A$, supported on the $\eta$-neighbourhood of $A$ and with $\Vert \nabla g \Vert_{L^{2}} \leq \eta$. 
\

{\bf Remarks}
\begin{itemize}\item  The terminology \lq\lq Minkowski measure'', \lq\lq capacity'' is provisional since the definitions may not exactly match up with standard ones in the literature. Since we only ever use the codimension 2 case we just talk about Minkowski measure  hereafter. 
\item In many  of our  applications the space $P$ will be a smooth Riemannian manifold outside $A$ and the meaning of   $\Vert \nabla  g\Vert_{L^{2}}$ is clear. In the general case this quantity can be defined by the theory in \cite{kn:CC2}\ \end{itemize}

\begin{prop}
If $m(A)$ is finite then $A$ has capacity zero. 
\end{prop}

To prove this, recall first that in this setting the volumes of $r$-balls in $P$ are bounded above and below by multiples of $r^{p}$. Thus the definition of Minkoski measure implies that the volume of the $r$-neighbourhood of $A$ is bounded by $C r^{2}$ for some fixed $C$.  For $\delta>0, Q>1 $ define a function $f_{\delta, Q}$ on $[0,\infty)$ by
\begin{itemize}
\item $f_{\delta, Q}(t)= 0$ if $0\leq t\leq \delta$;
\item $f_{\delta,Q}(t)=   \log (t/\delta)$ if $\delta\leq t\leq Q\delta$;
\item $f_{\delta, Q}(t)=\log Q $ if $t\geq Q\delta$.
\end{itemize}

Then $\vert f_{\delta, Q}'(t)\vert \leq \min( t^{-1}, \delta^{-1})$. Let $D_{A}: P\rightarrow [0,\infty)$ be the distance to $A$ and $f=f_{\delta, Q} \circ D_{A}$. So
$\vert \nabla f \vert \leq \min( D_{A}^{-1}, \delta^{-1})$. Let $S$ be the distribution function of $\vert \nabla f\vert^{2}$---so $S(\lambda)$ is the volume of the set  where $\vert \nabla f\vert^{2}\geq \lambda$. Thus $S(\lambda)$ is bounded by the volume of the $\lambda^{-1/2}$-neighbourhood of $A$ and so by $C \lambda^{-1}$. Also $S(\lambda)$ vanishes for $\lambda >\delta^{-2}$ and for all $\lambda$, $S(\lambda)$ is bounded by the volume of the $Q\delta$ neighbourhood of $A$, and so by $C Q^{2} \delta^{2}$. By the co-area formula
$$ \int_{P} \vert \nabla f\vert^{2} = \int_{0}^{\infty} S(\lambda)d\lambda\ =\int_{0}^{\delta^{-2}} S(\lambda) d\lambda. $$
Now $S(\lambda)\leq \min( C Q^{2} \delta^{2}, C\lambda^{-1})$ and we get
$$  \int_{P} \vert \nabla f \vert^{2} \leq C (1+ 2\log Q). $$
Now set $g= (\log Q)^{-1} (\log Q-f)$. So $g=1$ on the $\delta$-neighbourhood of $A$ and vanishes outside the $Q \delta$-neighbourhood. We have
$$\int_{P} \vert \nabla g \vert^{2} \leq C(\log Q)^{-2}(1+2\log Q). $$
Given $\eta>0$ we first make $Q$ so large that the right hand side in this formula is less than $\eta^{2}$ and then choose $\delta$ so small that $Q\delta<\eta$. This gives the desired function. 

We also record here:
\begin{prop}
Suppose that $A$ is closed and has Hausdorff dimension strictly less than $p-2 $. Then $A$ has capacity zero.
\end{prop}
This is proved in \cite{kn:DS}.

Beyond Proposition 1, another general  reason  for considering the notion of \lq\lq Minkowski measure''is that it behaves well with respect to Gromov-Hausdorff limits. Suppose we have metric spaces $X_{i}$ with Gromov-Hausdorff limit $Z$ and suppose given  any  subsets $A_{i}\subset X_{i}$. We fix distance functions $d_{i}$ on $X_{i}\sqcup Z$, as in the definition of Gromov-Hausdorff convergence. Then we  can define a \lq\lq limit'' $A_{\infty}\subset Z$ as the set of points $z\in Z$ such that there are $a_{i}\in A_{i}$ with $d_{i}(a_{i},z)\rightarrow 0$. However, without further information, we cannot say much about it. In particular the notion does not behave well with respect to Hausdorff measure.  For example we could  take $A_{i}$ to be finite and such that, even if we pass to any subsequence, $A_{\infty}$ is the whole of $Z$.  On the other hand, if we know that $A_{i}$ have bounded Minkowski measure, with a fixed bound, $m(A_{i})\leq M$ then it is straightforward to prove that (after taking a subsequence) the same is true in the limit. In  particular the limit has Hausdorff codimension at least 2.  

\subsection{Generalities on Gromov-Hausdorff limits}

As in Section 1, we consider a sequence of K\"ahler-Einstein metrics $\omega_{i}$ on $X_{i}$ with cone angle $2\pi\beta_{i}$ along $D_{i}$ where $\beta_{i}$ tends to $\beta_{\infty}<1$ and we suppose that these have a Gromov-Hausdorff limit $Z$. For convenience fix $\beta_{-}> 0$ and $\beta_{+}<1$ so that $\beta_{-}\leq\beta_{i}\leq \beta_{+}$ for all $i$.  We can apply the deeper results of Cheeger and Colding  to the structure of the limit space $Z$. Thus we have a notion of  tangent cones at  points of $Z$ (not {\it a priori} unique) and we can write
$Z=R\cup S$ where $R$ (the regular set) is defined to be the set of points where all tangent cones are $\bC^{n}$, and $S$ is defined to be $Z\setminus R$. We stratify the singular set as $$S=S_{1}\supset S_{2}\supset S_{3} \dots $$
  where $S_{j}$ is defined to be the set of points where no tangent cone splits off a factor $\bC^{n-j+1}$. The Hausdorff dimension of $S_{j}$ does not exceed $2(n-j)$. (Our labelling of the strata differs from what is standard in the literature, first because we use co-dimension rather than dimension and second because in this K\"ahler situation we can restrict to even dimensions.)
We can suppose that we have chosen fixed metrics $d_{i}$ on $X_{i}\sqcup Z$ realising the Gromov-Hausdorff convergence.  
\begin{prop}
The set $R\subset Z$ is open and the limiting metric induces a smooth K\"ahler-Einstein metric $\omega_{\infty}$ on $R$.
\end{prop}
Suppose that $q$ is a point in $D_{i}\subset X_{i}$. Then the limiting volume ratio (as the radius tends to $0$) at $q$ is clearly $\beta_{i}\leq \beta_{+}<1$. By Bishop-Gromov monotonicity the volume ratio of any ball centred at $q$ is less than $\beta_{+}$. Now let $p$ be a point of $R\subset Z$. By Colding's theorem on volume convergence it follows that there is some $\epsilon>0$ such that $d_{i}(q,p)\geq \epsilon$ for all such points $q$. Thus near $q$ we reduce to the standard theory of smooth K\"ahler-Einstein metrics.

\

\

The essential new feature of our case, as compared with the standard theory of smooth K\"ahler-Einstein metrics , is that $S\neq S_{2}$.
We  write $\cD= S\setminus S_{2}$ so by definition a point is in $\cD$ if it has a tangent cone which splits off $\bC^{n-1}$ but not $\bC^{n}$. Such a tangent cone must have the form $\bC_{\gamma}\times \bC^{n-1}$ where $0<\gamma<1$ and
 $\bC_{\gamma}$ denotes the standard flat cone of angle $2\pi{\gamma}$.
 \begin{prop}
 \begin{enumerate}
 \item Suppose a point $p\in {\cal D}$ has a tangent cone $\bC_{\gamma}\times \bC^{n-1}$ and a tangent cone $\bC_{\gamma'}\times \bC^{n-1}$. Then $\gamma=\gamma'$.
\item There is a $\gamma_{0}\in (0,1)$, depending only on the volume non-collapsing constant of $(X_{i},\omega_{i})$, such that if a point in ${\cal D}$ has tangent cone $\bC_{\gamma}\times \bC^{n-1}$ then $\gamma\geq \gamma_{0}$.
\end{enumerate}
\end{prop}
 
 These follow immediately from volume monotonicity, somewhat similar to Proposition 3. (In the second item recall that, as discussed in the Introduction, the metrics do satisfy a  fixed volume non-collapsing condition.)
 
 \
  
 All of the above discussion applies equally well to scaled limits. Thus, suppose we have a sequence $r_{i}\rightarrow \infty$ and base points $p_{i}\in X_{i}$. Then (taking a subsequence) we have a based Gromov-Hausdorff limit $Z'$ of
$(X_{i}, r_{i}^{2}\omega_{i}, p_{i})$ and the structure of the singular set of $Z'$ is just as described above. In particular we can apply this to tangent cones and iterated tangent cones of $Z$. Such a tangent cone  has the form $C(Y)$ where $Y$ is a $(2n-1)$-dimensional length space which has a smooth Sasaki-Einstein structure on an open set $R_{Y}$: the complement of $R$ has (real) Hausdorff codimension at least $2$ and can be written as a disjoint union ${\cal D}_{Y}\cup S_{2, Y}$ where $S_{2,Y}$ has (real) Hausdorff codimension at least $4$ and points of $\cD_{Y}$ have tangent cones $\bC_{\gamma}\times \bR \times \bC^{n-2}$.
 \
 \subsection{Review}

 In this subsection we assume knowledge of the main construction in \cite{kn:DS}, Section 3. Let $Z$ be a Gromov-Hausdorff limit of $(X_{i},\omega_{i})$, as above,  and let $p$ be a point in $Z$. Let $C(Y)$ be a tangent cone of $Z$ at $p$. 
\begin{defn}
We say that $C(Y)$ is a good tangent cone if the singular set $S_{Y}\subset Y$ has capacity zero.
\end{defn}

Now one of the main points in the argument of \cite{kn:DS} is
\begin{prop}
If ${\cal D}_{Y}$ is empty then $C(Y)$ is a good tangent cone.
\end{prop}
This follows from Proposition 2, since if ${\cal D}_{Y}$ is empty the singular set has Hausdorff codimension $\geq 4$. 
Of course by Proposition 1  we have:
\begin{prop}
If the singular set $S_{Y}$ has finite Minkowski measure then $C(Y)$ is a good tangent cone.
\end{prop}

Now we have:

\begin{prop}
Let $p$ be a point in a limit space $Z$ as above and suppose that there is a good tangent cone of $Z$ at $p$. Then there are $b(p), r(p)>0$ and an integer $k(p)$ with the following effect. Suppose $(X_{i}, D_{i},\omega_{i})$ is a sequence of K\"ahler-Einstein metrics with cone singularities and cone angles $\beta_{i}\rightarrow \beta_{\infty}<1$ then there is a $k\leq k(p)$ such that for sufficiently large $i$ there is a holomorphic section $s$ of $L^{k}\rightarrow X_{i}$, with $L^{2}$ norm $1$ and with $\vert s(x)\vert \geq b(p)$ for all points $x\in X_{i}$ with $d_{i}(x,p)<r(p)$.
\end{prop}

This statement is chosen to be exactly parallel to that of Theorem 3.2 in \cite{kn:DS} (except for trivial changes of language) and the proof is almost exactly the same. That is, the only facts from the convergence theory that  are needed for that proof are that the convergence is in $C^{\infty}$ on the regular set and that the singular set in $Y$  has capacity zero. The other point to note occurs in the $C^{0}$ and $C^{1}$ estimates for holomorphic sections (Proposition 2.1 in \cite{kn:DS}).
These can be obtained by approximating by metrics with positive Ricci curvature, or directly by working with the singular metrics and checking \lq\lq boundary terms''. The crucial point is that in the relevant Bochner-Weitzenbock formulae  the Ricci curvature enters with the favourable sign.

This gives the background to state our main technical result of this section:
\begin{thm}
With $Z$ as above, all tangent cones are good.
\end{thm}

Our proof of this will require  some  work and on the way we will establish a variety of other useful statements. 

\subsection{Gaussian sections}
We now want to refine the statement of Proposition 7 , which involves reviewing in more detail the constructions in \cite{kn:DS} and modifying them slightly.
When working with the line bundle $L^{k}\rightarrow X_{i}$ it is convenient to use the scaled metric $k\omega_{i}$. As in \cite{kn:DS} we use the notation
$L^{2,\sharp}$ etc. to denote quantities calculated using the rescaled metric. Likewise for the distance function $d^{\sharp}_{i}$ on $X_{i}\sqcup Z$.

\begin{prop}
Suppose that $p$ is a point in $Z$ with a good tangent cone. There is  a sequence $k_{\nu}\rightarrow \infty$ such that the following holds.  For any $\zeta>0$ there is an $R_{0}>1$, an integer $m\geq 1$ and integers $t_{i,\nu}$ with $1\leq t_{i,\nu}\leq m$ such that if $k= t_{i,\nu}k_{\nu}$ then for large enough $i$ there is a holomorphic section $s$ of $L^{k}\rightarrow X_{i}$ with following properties.
\begin{itemize} \item $ \vert \ \vert s(x)\vert- \exp(-d^{\sharp}_{i}(x,p)^{2}/4) \ \vert \leq \zeta$,  if $d^{\sharp}_{i}(x,p) \leq R_{0}$;
\item $\vert s(x)\vert \leq \zeta$ if  $d^{\sharp}_{i}(x,p) \geq R_{0}$;
\item $\Vert s \Vert_{L^{2,\sharp}}\leq (2\pi)^{n} + \zeta$.
\end{itemize}
 \end{prop}

\

We begin by recalling the basic strategy. First take a sequence of scalings $r_{\nu}\rightarrow \infty$ realising the given good tangent cone $C(Y)$ at $p$. Changing $r_{\nu}$ by bounded factors does not change the tangent cone so we can suppose that $r_{\nu}=\sqrt{k_{\nu}}$ for integers $k_{\nu}$.
For a point $z$ in the cone $C(Y)$ we use the notation $[z]$ to denote the distance to the vertex of the cone. There is a holomorphic section $\sigma_{0}$ of the trivial line bundle  over the smooth part of the cone, with $$\vert \sigma_{0}\vert= \exp(- [z]^{2}/4).$$
Thus the $L^{2}$ norm of $\sigma_{0}$ is $(2\pi)^{n/2}\sqrt{\kappa_{Y}}$ where $\kappa_{Y}\leq 1$ is the volume ratio of the cone. For $\eta'>0$, let $Y_{\eta'}$ be the complement of the $\eta'$-neighbourhood of the singular set in $Y$ and for $\delta<<1<<R$ let $U=U_{\eta', \delta, R}$ be the intersection of the cone on $Y_{\eta'}$ with the \lq\lq annulus'' defined by $\delta<[z]<R$. For suitable choice of parameters $\eta', \delta, R$ we construct an appropriate cut-off function $\chi$ of compact support on $U$ and equal to $1$ on a smaller set $U_{0}\subset U$. Then $\sigma=\chi \sigma_{0}$ is an \lq\lq approximately holomorphic section'' in that
\begin{equation}  \Vert \db \sigma \Vert_{L^{2}} \leq \Vert (\db \chi) \ e^{-[ z]\vert^{2}/4} \Vert_{L^{2}}, \end{equation}
and the right hand side will be chosen small. (In \cite{kn:DS} we denoted cut-off functions by $\beta$ but that symbol is taken here for the cone angle.)

We introduce some terminology. For $\epsilon>0$ we say that a map
$ \Gamma: U \rightarrow X_{i}$ is an {\it $\epsilon$-K\"ahler embedding} if it is a diffeomorphism onto  its image with the properties  that
\begin{equation}  \vert \Gamma^{*}(J_{i})- J\vert<\epsilon\ \ ,\ \ \vert \Gamma^{*}(k_{\nu}\omega_{i})-\omega\vert <\epsilon . \end{equation}
Here $J_{i}$ denotes the complex structure on $X_{i}$ and $J,\omega$ are the standard structures on the smooth part of the cone. By the definition of the tangent cone (and the smooth convergence on the regular set) we can, for any given $\epsilon$, choose $\nu$ sufficiently large and then $i\geq i(\nu)$ so that we have such a map. Further, for $r>0$ let $B_{\nu,i}(r)\subset X_{i}$ be the set of points $x_{i}\in X_{i}$ with distance $d^{\sharp}_{i}(x_i,p)\leq r$. Then we can define distance functions $d_{\nu,i}$ on $B_{\nu,i}(r)\sqcup \{ z: [z]\leq r\}\subset X_{i}\sqcup C(Y)$ realising the Gromov-Hausdorff convergence in the sense that for any fixed $r$ and for any $\delta>0$ then for sufficiently large $\nu$ and for $i\geq i(\nu)$ each component is $\delta$-dense.
Then (taking $r=2R$ say) we can also choose $\Gamma$ above such that 
\begin{equation} d_{i,\nu} (z, \Gamma(z)) \leq \epsilon. \end{equation}

Returning to our outline of the basic strategy, we choose such an $\epsilon$-K\"ahler embedding $\Gamma=\Gamma_{i}:U\rightarrow X_{i}$ for a suitable $\epsilon$.  A complication now arises  involving the potential \lq\lq topological obstruction'' from the holonomy of the connections. This handled in \cite{kn:DS} by the introduction of the factors $t_{i,\nu}$. We ignore this for the time being, so we suppose that $\Gamma$ is covered by a map of line bundles $\tilde{\Gamma}$ such that
\begin{equation} \vert \tilde{\Gamma}^{*}(A_{i})- A \vert <\epsilon \end{equation}
where $A$ is the preferred connection on the trivial bundle over the regular part of $C(Y)$ and $A_{i}$ is the connection on $L^{k}\rightarrow X_{i}$. 
With this done, we have a transplanted section $\tilde{\Gamma}_{*}(\sigma)$, which we  denote by $\sigma^{\sharp}$, of $L^{k}\rightarrow X_{i}$ with $L^{2,\sharp}$ norm approximately the same as that of $\sigma_{0}$ and with the $L^{2\sharp}$ norm of $\db \sigma^{\sharp}$ small. The basic \lq\lq H\"ormander construction'' is to pass to the holomorphic section $s=\sigma^{\sharp}-\tau$ with
$\tau = \db^{*}\Delta^{-1} \db \sigma^{\sharp}$, so that $\Vert \tau\Vert_{L^{2\sharp}}\leq \Vert \db \sigma^{\sharp}\Vert_{L^{2,\sharp}}$ will be small. Notice here that the third condition in  Proposition 8 is entirely straightforward. 

Next we recall that we have universal bounds on the $C^{0}$ and $C^{1}$ norms of holomorphic sections. We take these in local form; for a point $q$ in $X_{i}$ and a holomorphic section $s'$ of $L^{k}$ defined over a ball $B$ of fixed size in the re-scaled metric  (say the unit ball) centred at $q$ we have
\begin{equation} \vert s' (q) \vert \leq K_{0} \Vert s' \Vert_{L^{2, \sharp}(B)},\end{equation}
and
\begin{equation} \vert (\nabla s')(q) \vert \leq K_{1} \Vert s'\Vert_{L^{2, \sharp}(B)}. \end{equation}

The point to emphasise here is that $K_{0}, K_{1}$ are \lq\lq universal'' constants, which could be computed explicitly in terms of $\beta_{i}$ and the Hilbert polynomial. 

\

 With these preliminaries in place we can begin the proof of Proposition 8.  The argument is entirely elementary, but a little complicated. If written out in full the argument involves a large number of computable constants. For clarity,, we suppress these and often use a  notation such as \lq\lq  $<< K_{0}^{-1}$'', where strictly we mean \lq\lq $\leq$ some computable constant times $K_{0}^{-1}$''. We emphasise that while we have outlined the general scheme of construction above we are keeping the parameters $R, \eta', \delta, \epsilon$ free and they will be chosen below.

\

{\bf Step 1} Given $h>0$ there is some $\eta(h)\leq h$ such that for all $y$ in the singular set $S_{Y}$  there is a point $y'\in Y$ with $d(y,y')\leq h$ and such that $y'$ does not lie in the $2 \eta(h)$-neighbourhood of $S_{Y}$.

This follows from an elementary argument, using the facts that $S_{Y}$ is  closed and has empty interior.

\

{\bf Step 2} Choose $R_{0}>>1$ such that $\exp(-R_{0}^{2}/4)<<K_{0}^{-1} \zeta$. We will take $R\geq 2R_{0}$. 

\

{\bf Step 3} In the set up-as described above, with $U_{0}\subset U$ the set where $\chi=1$, it is clear that we have a fixed bound on $\vert \nabla \sigma^{\sharp}\vert $ over $\Gamma(U_{0})$. Using (8), and writing $\tau=s-\sigma^{\sharp}$ we get a bound of the form
\begin{equation} \vert \nabla  \tau\vert_{\chi(U_{0})}\leq K_{2}. \end{equation}
Now choose $\rho<< \zeta K_{2}^{-1}$.

\

{\bf Step 4}  Set $h=\rho/R_{0}$ and let $$\eta_{0}= \eta(h), $$
with $\eta(h)$ as in Step 1.

\

{\bf Step 5}
Now we recall that the cut-off function $\chi$ will have the form $\chi=\chi_{\eta',\eta}\ \chi_{\delta}\ \chi_{R}$ where:
\begin{itemize}
\item $\chi_{\eta',\eta}$ is a function on $Y$ (pulled back to the cone in the obvious sense), equal to $1$ outside the $\eta$-neighbourhood of the singular set and to $0$ on the $\eta'$-neighbourhood. Given $\eta$, use the definition of a good tangent cone to choose such a function so that the $L^{2}$ norm of its derivative is less than $\eta$ (and with $\eta'=\eta'(\eta)<\eta$, but otherwise uncontrolled).
\item $\chi_{\delta}$ is a standard cut-off function of $[ z]$, vanishing for $[ z] \leq \delta $ and equal to $1$ for $[z] \geq 2\delta$.
\item $\chi_{R}$ is a standard cut-off function of $[ z ]$, vanishing for $[ z] \geq R$ and equal to $1$ for $[ z] \leq R/2$.
\end{itemize}
This  fixes $\chi$, given parameters $\eta, \delta, R$: for the moment we have not fixed those parameters, but we  will take $\eta<\eta_{0}$.

\

{\bf Step 6}
Write $B(R_{0})$ for the set of $z\in C(Y)$ with $[ z] \leq R_{0}$ and $V(R_{0})$ for the intersection of the cone on $Y_{\eta_{0}}$ with the the annulus $2\delta < [ z] < R_{0}$. Thus, under our assumptions,
$\chi=1$ in $V(R_{0})$. Now, using Step 4, we can fix $\delta_{0}$ so that if $\delta\leq \delta_{0}$ then for all $z$ in $B(R_{0})$ we can find a $z'$ in $V(R_{0})$ such that:
\begin{itemize}
\item The distance between $z,z'$ does not exceed $\rho$;
\item $$ \vert \exp(-[ z]^{2}/4)-\exp(-[ z']^{2}/4) \vert << \zeta; $$
\item The $\eta_0\rho$-ball centred at $z'$ lies in $V(R_{0})$.
\end{itemize}

\

{\bf Step 7} Let $F_{0}$ be any function on $V(R_{0})$ with $|\nabla F_{0}\vert \leq  K_{2}$. (Here $F_{0}$ serves as a prototype for $\vert \tau\vert$.) If $z'\in V(R_{0})$ is any point as in Step 6 we get, by an elementary estimate, that there is a $K_{3}$ depending on the already chosen numbers $\zeta$, $\rho$, $\eta_0$ such that if $\Vert F_0 \Vert^{2}_{L^{2}} \leq K_{3}^{-1}$ then 
$ \vert F_0(z') \vert \leq \zeta/10$.

\

{\bf Step 8}
Let $F_{1}$ be any function on $B(R_{0})$ with $\vert \nabla F_{1}\vert \leq K_{1}$. (Here $F_{1}$ serves as a prototype for $\vert s\vert$).
Then for any point $z\in B(R_{0})$ we can find $z'$ as in Step 6 so we have
 \begin{equation}\vert F_1(z)-F_1(z')\vert << \zeta \end{equation}
 and also
 \begin{equation} \vert \exp(-[z]^{2}/4)- \exp(-[ z']^{2}/4)\vert << \zeta. \end{equation}

{\bf Step 9}
We now finally fix our parameters $\delta, \eta, R$. We choose $\delta<\delta_{0}, \eta<\eta_{0}, R>2 R_{0}$ but also we require that  
\begin{equation}   \Vert (\nabla\chi) e^{-[ z]^{2}/4}\Vert_{L^{2}} << \max\left( K_{3}^{-1},K_{0}^{-1} \zeta\right). \end{equation}

It is elementary calculation that this is possible (that is, the right hand side is a fixed number and we can make the left hand side as small as we please by taking $\eta $ and $ \delta$ small and $R$ large). 
 
 \

{\bf Step 10}
We are now in much the same position as in the corresponding stage of \cite{kn:DS}, with the difference that we have set up the background  to obtain more precise control of the final holomorphic section $s$. We choose an $\epsilon$-K\"ahler embedding $\Gamma$ of $U$ in $X_{i}$ satisfying (4), (5). Suppose first that we can find a lift $\tilde{\Gamma}$ as in (6). Then we define $\sigma^{\sharp}$ as described above. By making $\epsilon$ sufficiently small we can suppose that, following on from (12) we also have
\begin{equation}   \Vert \bar\partial \sigma^{\sharp}\Vert_{L^{2,\sharp}} << \max\left( K_{3}^{-1},K_{0}^{-1} \zeta\right). \end{equation}
 At this point we can fix the parameter $m$. This goes just as in \cite{kn:DS}: we need to arrange that after dilating by a suitable factor $t$ with $1\leq t\leq m$ the potential topological obstruction arising from holonomy is a small perturbation, so we can construct a section satisfying the constraints (13). To simplify exposition,  assume that in fact $t=1$.

\

{\bf Step 11}  For suitably large $\nu$ and $i\geq i(\nu)$ we can suppose that for any point $x\in X_{i}$ with $d^{\sharp}(x,p)\leq R_{0}$ there is a point $x'$ in   $\Gamma_{i}(V)$ with $d^{\sharp}(x,x')\leq 2\rho$ and such that the ball of radius $\eta_0\rho$ centred at $x'$ lies in $\Gamma_{i}(V)$.
Here $V=V(R_{0})$. Now we apply the estimates, modelled on the prototypes in (10),(11), to deduce that $\vert \vert s(x)\vert - \exp(d^{\sharp}_{i}(x,p)^{2})\vert \leq \zeta$.

\

{\bf Step 12}
The remaining task is to obtain the estimate $\vert s(x) \vert\leq \zeta$ when $d^{\sharp}_{i}(x,p) \geq R_{0}$. For this we use the $C^{0}$ estimate (7). If $x\in X_{i}$ is such a point  then the $L^{2,\sharp}$ norm of $\sigma^{\sharp}$ over the unit ball centred at $x$ is small (compared with $\zeta$) by Step 1. Likewise for the $L^{2,\sharp}$ norm of $\tau$ by Step 9. Then we obtain the conclusion from the $C^{0}$ estimate.

\

{\bf Note.} In our applications below we can take $\zeta$ to be some fixed small number, say $1/100$. Then, given $p$, the number $m$ is fixed. Using a diagonal argument and passing to a subsequence we may suppose that all the $t_{\nu,i}$ are equal and this effectively means that we can ignore this extra complication.

\

\

If we have suitable holomorphic functions  on the cone we can construct more holomorphic sections. 
\begin{prop}Suppose that $p$ is a point in $Z$ with a good tangent cone and that $f$ is a holomorphic function on the regular part of the cone such that for some $\alpha\geq 1$ and $C>0$ we have: \begin{itemize} \item  There are smooth functions $G_{\pm}$ of one positive real variable with $\vert G_{\pm}(r)\vert \leq C r^{\alpha} , \vert G'_{\pm}(r)\vert \leq C r^{\alpha-1}$ and 
\begin{equation}   G_{-}([z])\leq \vert f (z)\vert \leq G_{+}([z]). \end{equation}
\item $\vert \nabla f (z) \vert \leq C [z]^{\alpha-1}$.
\end{itemize}
 
Let $k_{\nu}$ be the sequence as in Proposition 8.  For any $\zeta>0$ we can choose  $R_{0}, m, t_{i,\nu}$ as in Proposition 8 so that in addition  there is a holomorphic section $s_{f}$ of $L^{k}\rightarrow X_{i}$ such that, writing $d= d^{\sharp}_{i}(x,p)$
\begin{itemize} \item $ \exp(-d^{2}/4)G_{-}(d)-\zeta\leq  \vert s_{f}(x)\vert \leq \exp(-d^{2}/4) G_{+}(d)+\zeta$  if $d \leq R_{0}$;
\item $\vert s_{f}(x)\vert \leq \zeta$ if  $d \geq R_{0}$;
\item $\Vert s_{f} \Vert_{L^{2,\sharp}}\leq N +\zeta$';
\end{itemize}
where $N$ can be computed explicitly from $C, \alpha$. 

\end{prop}
While the statement is a little complicated the proof is essentially identical to that of Proposition 8, starting with the section $\chi f \sigma_{0}$. The Gaussian decay of $\sigma_{0}$ dominates the polynomial growth of $f$ and the condition $\alpha\geq 1$ means that there are no problems at the vertex of the cone. 

For $r>0$ we write $B^{\sharp}=B^{\sharp}(p,r)\subset X_{i}$ for the set of points  $x$ with $d^{\sharp}(x,p)< r$.  (For our purposes below one could take  $r=3$,   say.)
Given any fixed $r$ then clearly if $\zeta$ is chosen sufficiently small the radius $R_{0}$ in Proposition 8 will be much bigger than $r$ and we will have a definite lower bound $\vert s\vert \geq c>0$ for our section over $B^{\sharp}$. If then we have a function $f$ as in Proposition 9, we can form the quotient $\tilde{f}= s_{f}/s$ as a holomorphic function over $B^{\sharp}$. This satisfies a fixed $C^{1}$ bound
\begin{equation}
 \vert f \vert, \vert \nabla f \vert \leq M(C,\alpha,K_{0}, K_{1}, c)
  \end{equation} 
for a computable function $M$.

We can say more.\begin{prop}
 Let $W$ be any pre-compact open subset of the regular set in the ball $\{[z]<r\}$ in $C(Y)$. For any $\epsilon, \zeta>0$, for $\nu$ sufficiently large and for $i\geq i(\nu)$ there is an $\epsilon$-K\"ahler embedding $\Gamma:W\rightarrow X_{i}$ such that $\vert f - \tilde{f}\circ \Gamma\vert\leq \zeta$ on $W$.
\end{prop}
 
 The proof is again essentially the same as Proposition 8.

\

 Now we  consider the extension of all these ideas to scaled limits. So suppose that we have a sequence $r_{i}\rightarrow \infty$, base points $q_{i}\in X_{i}$ and that $Z'$ is the based limit of $(X_{i}, r_{i}^{2}\omega_{i})$.
There is no loss of generality in supposing $r_{i}^{2}= a_{i}$, for integers $a_{i}$. Thus, algebro-geometrically, the scaling corresponds to considering line bundles $L^{a_{i}}\rightarrow X_{i}$. Let $d_{i}$ now denote a distance function for the re-scaled metrics, realising the Gromov-Hausdorff convergence.
Given $k$ we write, as before, $d_{i}^{\sharp}= k^{1/2} d_{i}$, and in general as before we add a $\sharp$ to denote quantities calculated in the metric scaled by $\sqrt{k}$. The point is that this is {\it additional} to the scalings we have already made by $\sqrt{a_{i}}$.
\begin{prop}
Suppose that $p$ is a point in $Z'$ with a good tangent cone. There is  a sequence $k_{\nu}\rightarrow \infty$ such that the following holds.  For any $\zeta>0$ there is an $R_{0}>1$, an integer $m\geq 1$ and integers $t_{i,\nu}$ with $1\leq t_{i,\nu}\leq m$ such that if $k= t_{i,\nu}k_{\nu}$ then for large enough $i$ there is a holomorphic section $s$ of $L^{a_{i} k}\rightarrow X_{i}$ such that
\begin{itemize} \item $ \vert \ \vert s(x)\vert- \exp(-d^{\sharp}_{i}(x,p)^{2}/4) \ \vert \leq \zeta$  if $d^{\sharp}_{i}(x,p) \leq R$;
\item $\vert s(x)\vert \leq \zeta$ if  $d^{\sharp}_{i}(x,p) \geq R$;
\item $\Vert s \Vert_{L^{2,\sharp}}\leq (2\pi)^{n} + \zeta$.
\end{itemize}
 \end{prop}
 This is {\it precisely} the same statement as Proposition 8 except that we have $Z'$ in place of $Z$ and $L^{a_{i} k}$ in place of $L^{k}$. (So we can regard Proposition 8 as a special case when all $a_{i}=1$.) The point is that the proof is {\it precisely} the same. (It may be pedantic to have written out the whole statement but that seems the best way to be clear.) Likewise  the statements in Propositions 9,10 apply  to scaled limits, introducing the extra powers $a_{i}$.

\subsection{Points in ${\cal D}$}

Next we focus attention on a point $p\in {\cal D}\subset Z$. By definition this means that there is some tangent cone $\bC_{\gamma}\times \bC^{n-1}$.  This is clearly a good tangent cone so the results above apply. We take complex co-ordinates $(u, v_{1}, \dots, v_{n-1})$ on $\bC_{\gamma}\times \bC^{n-1}$ with the metric 
\begin{equation}  \vert u\vert^{2\gamma-2}\vert du\vert^{2} + \sum \vert dv_{i}\vert^{2}. \end{equation}
Then if $z=(u,v_{1}, \dots, v_{n-1})$ we have, in our notation above, 
\begin{equation} [z]^{2} =|\gamma|^{-2}\vert u \vert^{2\gamma}+ \sum \vert v_{i}\vert^{2}. \end{equation}
Each of the co-ordinate functions $u, v_{i}$ satisfy the hypotheses of Proposition 9 so we get  holomorphic functions $\tilde{u}, \tilde{v}_{i}$ on $B^{\sharp}$. Now we regard these as the components of a map $$F:B^{\sharp}\rightarrow \bC^{n}. $$
An easy extension of Proposition 9 shows that for any $\zeta>0$ we can suppose that 
\begin{equation} \vert [F(x)]- d(x,p) \vert \leq \zeta, \end{equation}
where we regard $[\ ]$ as the function on $\bC^{n}$ defined by (17). 
 Proposition 10 implies that for any precompact subset $W$ in the set
 $$\{z=(u ,v_{1}, \dots, v_{n-1}): [z]<r, u\neq 0\}$$
 and any $\epsilon, \zeta>0$ we can suppose that there is an $\epsilon$-K\"ahler embedding $\Gamma:W\rightarrow  B^{\sharp}$ such that
$$  \vert F\circ \Gamma(z) - z\vert \leq \zeta, $$
for all $z$ in $W$. (Here, of course, the meaning of \lq\lq we can suppose that'' is that we should take $\nu$ sufficiently large and $i\geq i(\nu)$.)

Now fix $r=3$. For points $z\in \bC^{n}$ with $\vert z\vert^{2}\leq 1$ we have an elementary inequality $[z]<2$ so if $\zeta<1/10$ we have $\vert f(x)\vert >1$ on the boundary of $B^{\sharp}$. Let $\Omega=\Omega_{i}\subset X_{i}$ be the preimage of the open unit ball in $\bC^{n}$:
$$  \Omega=\{ x\in B^{\sharp}\subset X_{i}: \vert F(x)\vert < 1\}. $$ 
 
 The next result is one of the central points in this paper. Later, we will develop it further to show that $Z$ has an natural complex manifold structure around points of ${\cal D}$.  

\begin{prop} For $\nu$ sufficiently large and $i\geq i(\nu)$ the map $F$ is a holomorphic equivalence  from $ \Omega_{i}\subset X_{i}$ to the unit ball $B^{2n}\subset \bC^{n}$. 
\end{prop}
 By construction $F:\Omega\rightarrow B^{2n}$ is a proper map and so has a well-defined degree. The fibres of $F$ are compact analytic subsets but we have constructed a non-vanishing section $s$ of $L^{k}$ over $\Omega\subset B^{\sharp}$ so the fibres must be finite. To establish the proposition it suffices to show that the degree is $1$. By $(15)$ we  have a fixed bound on the  derivative of $F$ over $\Omega$ , say $\vert \nabla F \vert \leq K_{4}$.

Recall that we are making the identification $\bC^{n}= \bC_{\gamma}\times \bC^{n-1}$. We write $d_{\gamma}(z,z')$ for the distance between points $z,z'$ in the singular metric, not to be confused with Euclidean distance $\vert z-z'\vert$. 

Let $z_{0}$ be the point $(1/2,0,\dots, 0)\in B^{2n}\subset \bC^{n}$. The distance, in the singular metric, from $z_{0}$ to the boundary of $B^{2n}$ is
$\gamma^{-1}(1- (1/2)^{\gamma})= d$, say. The distance in the singular metric from $z_{0}$ to the singular set $S=\{u=0\}$ is $\gamma^{-1}(1/2)^{\gamma}>d$. Fix a small number $\rho>0$ and a pre-compact open set $W\subset B^{2n}\setminus S$ with the following properties.
\begin{enumerate}
\item $W$ contains the Euclidean ball $B(z_{0}, 2\rho)$ with centre $z_{0}$ and radius $2\rho$.
\item $d_{\gamma}(B^{2n}\setminus W, B(z_{0}, 2\rho))> d/2$;
\item For any $z\in B^{2n}$ there is a $z'$ in $W$ with $d_{\gamma}(z,z')< \delta_{0}= \min( \rho/(20 K_{4}), d/4)$.
\end{enumerate}
Now choose $\epsilon, \zeta$ such that 
\begin{itemize}
\item $\zeta<\rho/2$;
\item $K_{4}(\delta_0 + 2\epsilon)+ \zeta< \rho/10$;
\item $ \delta_0< d/2-4\epsilon$.
\end{itemize}

 We can suppose that there is an $\epsilon$-K\"ahler embedding $\Gamma:W\rightarrow X_{i}$ such that
$\vert F\circ \Gamma(z)- z\vert \leq \zeta<\rho/2$ for $z\in W$. 
 This implies that $F\circ \Gamma$ maps the boundary of $B(z_{0}, \rho)$ to $\bC^{n}\setminus \{z_{0}\}$. This boundary map has a well defined degree (defined by the action on $(2n-1)$-dimensional homology) and the degree is $1$. Now once $\epsilon$ is reasonably small it is clear from the definition that the derivative of $\Gamma$ is invertible at each point and orientation preserving. Then it follows by basic differential topology and complex analysis that for any point $z$ with $\vert z- z_{0}\vert \leq \rho/4$ there is a unique $\tilde{z}$ with $\vert \tilde{z}-z_{0}\vert \leq \rho$ such that $F\circ\Gamma(\tilde{z})=z$.
In particular we could take $z=z_{0}$, so we get a $\tilde{z}_{0}$ with $\vert \tilde{z}_{0}-z_{0}\vert \leq \rho$ and $F\circ \Gamma(\tilde{z}_{0})= z_{0}$.
What we have to show is that the only point $x\in \Omega$ with $F(x)= z_{0}$ is $ \Gamma(\tilde{z}_{0})$. 

Suppose then that $x\in \Omega$ and $F(x)= z_{0}$. By (5) and item 3 above
we can find an $x'$ in $\Gamma(W)$ with $d^{\sharp}(x,x')\leq \delta_{0}+ 2\epsilon$. Thus $\vert F(x)-F(x')\vert \leq K_{4}(\delta_{0}+2\epsilon)$.
Writing $x'=\Gamma(z')$ we have $\vert z_{0}- z'\vert \leq \vert z_{0}- F\circ \Gamma(z')\vert + \vert z'-F\circ \Gamma(z')\vert \leq K_{4}(\delta+2\epsilon)+ \zeta < \rho/10$.  Thus $z'$ lies in $B(z_{0}, \rho)$ and so $d_{\gamma}(z', B^{2n}\setminus W)\geq d/2$ and by (5) again  $d^{\sharp}(x', X_{i}\setminus \Gamma(W))\geq d/2-2\epsilon$. Since $\delta_{0}< d/2-4\epsilon$ we see that $x$ itself lies in $\Gamma(W)$ so we can take $x=x'$ above and we have $x=\Gamma(z')$ with $F\circ\Gamma(z')= z_{0}$. Thus 
$\vert z'-z_{0}\vert\leq\zeta<\rho/2$ and by the uniqueness statement in the previous paragraph we must have $z'=\tilde{z}_{0}$.

\

\

We can use the map $F$ to regard $k\omega_{i}$ as a metric $\tomega_{i}$ on the unit ball $B^{2n}\subset \bC^{n}$, with co-ordinates $(u, v_{1}, \dots, v_{n-1})$.  We also write $\tD_{i}$ for $F_{i}(\Omega_{i}\cap D_{i})$.
Let $N_{\delta}$ be the (Euclidean) $\delta$-neighbourhood of the hyperplane
$\{ u=0\}$. Let $\omega_{\Euc}$ be the Euclidean metric and $\omega_{(\gamma)}$ be the standard cone metric representing $\bC_{\gamma}\times \bC^{n-1}$. To sum up what we know: by making $\nu$ sufficiently large and then $i\geq i(\nu)$  we can arrange the following, for fixed $C$.
\begin{itemize}
\item $\tomega_{i}=\dbd \phi_{i}$ where $0\leq \phi_{i}\leq C$;
\item $\tomega_{i}\geq C^{-1} \omega_{{\Euc}}$;
\item For any $\zeta,\delta$ we can suppose $\vert \tomega_{i}-\omega_{(\gamma)} \vert \leq \zeta $ outside $N_{\delta}$ (and likewise for any given number of derivatives). 

\end{itemize}

\subsection{Consequences}

We will work in the co-ordinates $(u,v_{1}, \dots, v_{n-1})$ on the unit ball  in $\bC^{n}$, as above. It is clear that the singular sets $\tD_{i}$ must lie in the tubular neighbourhood $N_{\delta}$. This means that there  are well-defined homological intersection numbers $\mu=\mu_{i}$ with  a disc transverse to the hyperplane $\{u=0\}$. The singular set is given, in this co-ordinate chart by a Weierstrasse polynomial
\begin{equation} u^{\mu} + \sum_{j=0}^{\mu-1} A_{j} u^{j}=0, \end{equation}
where the $A_{j}$ are holomorphic functions of $v_{1}, \dots , v_{n-1}$. (Since the precise size of this domain is arbitrary we can always suppose that picture extends to a slightly larger region.)
\begin{prop}
We have the identity $(1-\gamma)=\mu_{i}(1-\beta_{\infty})$. \end{prop}

A first consequence of this is that $\mu$ does not depend on $i$. A second consequence is that there are only a finite number of possibilities for $\gamma$;
that is  $\gamma=  1-\mu(1-\beta_{\infty})$ for $\mu=1,\dots ,\mu_{{\rm max}}$ where $\mu_{{\rm max}}$ is the largest integer $\mu$ such that $1-\mu(1-\beta_{\infty})>0$.
So we can write \begin{equation}{\cal D}=\bigcup_{\mu=1}^{\mu_{\rm max}} {\cal D}_{\mu}.\end{equation}

\

{\bf Discussion}
The possible existence of points in ${\cal D}_{\mu}$ for $\mu>1$ is the major difficulty in the limit theory for K\"ahler-Einstein metrics with cone singularities.
As we shall see in (3.2), we understand the points in ${\cal D}_{1}$ essentially as well as the regular points. But one expects that there will be cases when the ${\cal D}_{\mu}$ for $\mu>1$ arise. For $\epsilon>0$ and $\beta>1/2$, cut out two disjoint wedge-shaped regions from $\bC$ with angles $2\pi(1-\beta)$ and vertices at the points $\pm \epsilon$. Then identify the boundaries of the resulting space to get a manifold with two cone singularities of angle $2\pi\beta$. The limit as $\epsilon$ tends to $0$ has a single singularity with cone angle $2\pi( 2\beta-1)$. This is the model for the simplest way in which points of ${\cal D}_{2}$ can appear, but one expects that more complicated phenomena  arise.  Similarly, one should probably  be careful in understanding the structure of the sets ${\cal D}_{\mu}$. Let $D$ be the cusp curve
with equation $u^{2}=v^{3}$ in $\bC^{2}$. Then it seems likely that for $\beta>1/2$ there is a Ricci flat K\"ahler metric on the complement of the curve with standard cone singularities of cone angle $2\pi\beta$ away from the origin, so the tangent cones at these points would be $\bC_{\beta}\times \bC$. It seems then likely (or, at least, hard to rule out) that the tangent cone at the origin is $\bC_{\gamma}\times \bC$ with $\gamma=(2\beta-1)$, i.e. an isolated point of ${\cal D}_{2}$.  

\

To prove Proposition 13, we consider a disc $T$ defined by fixing $v_{j}$ at some fixed values close to $0$ and constraining $\vert u\vert \leq 1/2 $ say.
We can suppose that the disc meets $\tD_{i}$ transversally in $\mu$ points.  Let $\rho$ be the Ricci form of $\tomega_{i}$. 
 \begin{lem} Let ${\rm Hol}_{i}\in S^{1} \subset \bC$ denote the holonomy of the anti-canonical bundle around the boundary of $T$. Then 
\begin{equation}   {\rm Hol}_{i} = \exp\left( 2\pi \sqrt{-1}( \mu (1-\beta_{i}) + \int_{T} \rho ) \right) . \end{equation}
\end{lem}

Choose local complex co-ordinates $z_{1}, \dots, z_{n}$ in a small neighbourhood of a point of $T\cap D_{i}$ so that $D_{i}$ is defined locally by the equation $z_{1}=0$ and $D_{i}$ by  $z_{2}= \dots z_{n}=0$. Let ${\rm Hol}(r)$ denote the holonomy of the anticanonical line bundle around the circle $\vert z_{1}\vert^{\beta}=r$ in $T$. By the Gauss-Bonnet formula it suffices to prove that, for all such intersection points, ${\rm Hol}(r)$ tends to $\exp(1-\beta_{i})$ as $r$ tends to zero.  Let $L$ be the logarithm of the ratio of the volume form of $\omega$ and $\vert z_{1}\vert^{-2\beta}$ times the Euclidean volume element, in these co-ordinates. The hypothesis that $\omega_{i}$ has a standard cone singularity implies that $L$ is a bounded function. Take polar co-ordinates $r,\theta$ on $D$, with $r=\vert z_{1}\vert^{\beta}$ and $\theta$ the argument of $z_{1}$. Then
$$   {\rm Hol}(r)= \exp(1-\beta_{i} \exp(I(r))$$
where $$ I(r)= r \int \frac{\partial L}{\partial r} d\theta, $$
The fact that the Ricci curvature of $\omega_{i}$ is bounded implies, using the Gauss-Bonnet formula again, that
$$I(r_{1})-I(r_{2})\vert \leq C (r_{1}^{2}- r_{2}^{2}), $$
so $I(r)$ has a limit $h$ as $r$ tends to $0$.  Thus
$$    \frac{d}{d r}\int L d\theta = h r^{-1} + o(r^{-1}), $$
and if $h$ is not zero this implies that $\int L d\theta$ is unbounded as $r\rightarrow 0$ which  is a contradiction.

 Note that the point in the proof of the Lemma is that we can work for a fixed $i$ on an arbitrarily small neighbourhood of an intersection point. We are not claiming that there is any uniformity in the size of this neighbourhood as $i$ varies. 

\

 Making $\nu$ and $i$ large we can make ${\rm Hol}_{i}$ as close as we please to the corresponding holonomy for the metric $\omega_{(\gamma)}$, which is
$\exp(2\pi \sqrt{-1}(1-\gamma))$. On the other hand $\rho=k^{-1} (1-\lambda(1-\beta_{i})) \tomega_{i}$ where $k$ can be made as large as we want by making $\nu $ large. So (using  Lemma 1) to establish Proposition 13 it suffices to prove that the area of $T$, in the  metric $\tomega_{i}$, satisfies a fixed bound. 

Now we use a version of the  the well-known Chern-Levine-Nirenberg inequality.
We fix a standard cut-off function $G$ of $\vert u\vert$, equal to $1$ when
$\vert u\vert\leq 1/2$ and vanishing for $\vert u\vert\geq 3/4$. Then
$$  \int_{T} \tomega_{i} \leq \int_{2T} G \omega^{\sharp}_{i}= \int_{2T} \dbd G \phi_{i} \leq C \int_{2T} \vert\dbd G\vert $$
where $2T$ has the obvious meaning and we use the bound on the K\"ahler potential $\phi_{i}$.   

{\bf Remarks}
\begin{enumerate}
\item It may be possible to prove Proposition 13  by adapting the Cheeger-Colding \lq\lq slicing'' theory using harmonic functions. But we avoid that by using the \lq\lq holomorphic slicing'' as above.
\item A similar argument, applied to the Chern connection on the line bundle $L^{k}$, shows that in fact the potential \lq\lq topological obstruction'' arising from holonomy does not occur. So for points of ${\cal D}$ we do not need to introduce the scale factors $t_{i,\nu}$.
\end{enumerate}

In the next two propositions we study the volume of the divisor $\tD_{i}\subset B^{2n}$, with respect to the metric $\tomega_{i}$

\begin{prop} There is a universal constant $c_{1}>0$ such that
$$  \int_{\tD_{i}} \tomega_{i}^{2(n-1)} \geq c_{1}.$$
\end{prop}

Since the Euclidean metric is dominated by a multiple of $\tomega_{i}$ it suffices to get a lower bound on the Euclidean volume. But this is clear. Recall that if $A$ is any $p$-dimensional complex analytic variety in $B^{2n}\subset \bC^{n}$ and if we write $V(r)$ for the volume of the intersection of $A$ with the $r$ ball then $ r^{-2p}V(r)$ is an increasing function of $r$. If the origin lies in $A$ then the limit of this ratio as $r$ tends to $0$ is an least the volume of the ball in $B^{2p}$, so we get a lower bound on $V(1)$. In our case the divisor $\tD_{i}$ need not pass exactly through the origin but it contain a point within the $\delta$ neighborhood of the origin (in the Euclidean metric) so the same argument applies.

Next let $\frac{1}{2} B^{2n}$ denote as usual the ball of radius $1/2$. 
\begin{prop}
There is a universal constant $c_{2}$ (depending on $\mu$) so that
$$  \int_{\frac{1}{2} B^{2n}\cap \tD_{i}} \tomega_{i}^{2(n-1)} \leq c_{2}. $$
\end{prop}
(Of course, by adjusting our set-up slightly, we could replace $\frac{1}{2} B^{2n}$ by any given sub-domain.)

\

The proof  uses another variant of the Chern-Levine-Nirenberg argument. We fix a sequence of cut-off functions $G_{1}, \dots, G_{n-1}$ depending on $\vert \uv\vert$ (where $\uv =(v_{1}, \dots, v_{n-1})$), such that $G_{n-1}$ vanishes when $\vert \uv \vert \geq 3/4$, $G_{1}=1 $ when $\vert v\vert \leq 1/2$ and $G_{j+1}=1$ on the support of $G_{j}$, for $j=1,\dots, n-2$.
Let 
$$ \Theta= \sqrt{-1} \sum_{j=1}^{n-1} dv_{i} d\overline{v}_{i} $$
i.e the Euclidean form pulled up from the $\bC^{n-1}$ factor. Then we can suppose that $\dbd G_{j} \leq C \Theta$ for some fixed $C$ (since the $G_{j}$ depend only on $\uv$). On the other hand it follows from our set-up that, if $\delta$ is reasonably small, the restrictions of $G_{j}$ to $D_{i}$ are functions of compact support. Now we write
$$  \int_{\frac{1}{2} B^{2n}\cap \tD_{i}} \tomega_{i}^{n-1} \leq  \int_{ B^{2n}\cap \tD_{i}} G_{1} (\dbd\phi_{i})^{n-1}, $$
Integrating by parts,  the right hand side is bounded by
$$  C  \int_{B^{2n}\cap \tD_{i}} \phi_{i} \Theta (\dbd \phi_{i})^{n-2}, $$
and we continue in the usual fashion to exchange factors of $\dbd\phi_{i}$ for factors of $\Theta$. After $(n-1)$ steps we get a bound in terms of
$$  \int_{\tD_{i}\cap \{ \vert \uv \vert \leq 3/4\}} \Theta^{n-1}, $$
which is exactly $\mu$ times the volume of the $3/4$-ball in $\bC^{n-1}$.

All of the discussion above is local and applies just as well to scaled limits
$Z'$, in the manner described in (2.5) above. 
\subsection{Lower bound on densities}

Let $\vert \ \vert_{*}$ denote the norm on $\bC^{q}$
$$  \vert (y_{1}, \dots, y_{q}) \vert_{*}= \max_{j} \vert y_{j}\vert. $$

\begin{prop} Suppose that $p$ is a point in a scaled limit space $Z'$ which is a based limit of $(X_{i}, a_{i} \omega_{i})$ for integers $a_{i}$. Let $d_{i}$ be distance functions as usual. Suppose  that there is a neighbourhood $N\subset Z' $ of $p$ such that for all $p'$ in $N$ there is at least one good tangent cone to $Z'$ at $p'$.  Then we can find  $C$ and $ \rho_{1}> \rho_{2}>0$ and $q$ such that the following is true. For large enough $i$ there is a holomorphic map
$F:\Omega_{i}\rightarrow \bC^{q}$ where $\Omega_{i}\subset X_{i}$ is an open set containing the points $x\in X_{i}$ with $d_{i}(p,x)\leq \rho_{1}$ and
\begin{itemize}
\item $\vert \nabla F_{i}\vert \leq C$ (where the norm is computed using the scaled metric $a_{i}\omega_{i}$);
\item For points $x$ with $d_{i}(p,x)=\rho_{1}$ we have $\vert F_{i}(x)\vert_{*}\geq 1/2$;
\item For points $x$ with $d_{i}(p,x)\leq \rho_{2}$ we have $\vert F_{i}(x)\vert_{*}\leq 1/100$.
\end{itemize}
\end{prop}

We apply Proposition 11 to find  $k, R$ and a holomorphic section $s$ of $L^{k a_{i}}\rightarrow X_{i}$ (for large $i$) such that $\vert \vert s(x)\vert - \exp(-d^{2}/4)\vert \leq 1/100$ if $d\leq R$ and $\vert s(x) \vert \leq 1/100$ if $d\geq R$. Here $d=\sqrt{k}d_{i}(p,x)$. We can take $k$ arbitrarily large so there is no loss in supposing that all points in $Z'$ a distance less than $ 2 k^{-1/2}$ from $p$ have good tangent cones. There is also no loss of generality (multiplying by a factor slightly larger than $1$ if necessary)  in supposing that there is some fixed $t_{0}>0$ such that $\vert s\vert \geq 1$ at points where $d\leq t_{0}$. 

Let $t_{+}, t_{-}$ be the fixed numbers such that
$$  \exp(-t_{\pm}^{2}/4)=1/2\mp 1/100,$$
so $t_{+}> t_{-}$.
Let ${\cal A}\subset Z'$ be the closed annulus of points distance between $t_{-}, t_{+}$ from $p$. Since $t_{\pm}< 2$, each point $p'$ of ${\cal A}$ has some good tangent cone. Applying Proposition 11 we can find some integer $l(p')\geq 10 $ and a number $R(p')$ such that for all sufficiently large $i$ there is a  holomorphic section $\sigma_{p'}$ of $L^{ k l(p') a_{i}}$ which obeys just the same estimates as above but with $d=\sqrt{k l(p')} d_{i}$. (That is, having fixed $k$ we can replace $L$ by $L^{k}$ and then apply Proposition 11.) By compactness of ${\cal A}$ we can find a finite collection $p'_{j}$ for $j=1,\dots, q$ such that the balls of radii $r_{j}= \frac{t_{+}}{2} \sqrt{k l(p'_{j})}$ with these centres cover ${\cal A}$. Now define $F_{i}$, mapping into $\bC^{q}$, to be the function
with components $\sigma_{p'_{j}}/ s^{l(p'_{j})}$. Let $\rho_{1}= k^{-1/2} t_{+}$. Then for points $x$ with $d_{i}(x,p)\leq \rho_{1}$ we have $\vert s(x)\vert\geq 1/2$ so $F_{i}$ is well-defined at such points and satisfies a $C^{1}$ bound. For sufficiently large $i$, for any point $x$ with $d_{i}(x,p)=\rho_{1}$ there is some $j$ such that $d_{i}(x,p'_{j})\leq \frac{t_{+}}{2} \sqrt{k l(p'_{j})}$ and this means that $\vert \sigma_{p'_{j}}\vert \geq 1/2$. Thus for such a point
$$   \vert F_{i}(x) \vert_{*} \geq 1/2 , $$
since it is clear that $\vert s \vert\leq 1$ at such points. Let $\rho_{2} = k^{-1/2} t_{0}$ so at points $x$ with $d(x,p)\leq \rho_{2}$ we have $\vert s(x)\vert \geq 1$. By our choice that $l_{j}\geq 10$ it is clear that $\vert \sigma_{p'_{j}}(x)\vert\leq 1/100$ for each $j$. Thus at such points
$$ \vert F_{i}(x) \vert_{*} \leq 1/100. $$

 For a point $x\in X_{i}$ and $r>0$ define the density function 
 \begin{equation}  V(x,r)= r^{2-2n} {\rm Vol}( D_{i} \cap B_{r}(x)), \end{equation}
 where of course the volume is computed using $\omega_{i}$.

 \begin{prop} There is a constant $c>0$ such that for all $i$, all points $x\in D_{i}\subset X_{i}$ and all $r\leq 1$ we have $V(x,r)\geq c$.
\end{prop}
The proof is by contradiction, so suppose we have  sequences $x_{i}\in D_{i}$ and $r_{i}\leq 1 $ such that $V(x_{i}, r_{i})\rightarrow 0$. Rescale by factors
$r_{i}^{-1}$ and take the based limit space $(Z',p)$. We claim first that every tangent cone to $Z'$ at a point $z$ in $Z'$ of distance strictly less than $1$ from $p$ is good. In fact we claim that if $C(Y)$ is such a tangent cone at $z$ then ${\cal D}(Y)$ is empty. Given this the assertion that the tangent cone is good follows from Proposition 5. 

Suppose that there is some $y\in {\cal D}_{Y}$. The tangent cone $C(Y)$ is itself a scaled limit of the $X_{i}$, for a suitable choice of scalings and base points. Thus we can apply Proposition 14. It follows that there is some definite lower bound on the volume of the singular set $D_{i}$ in the unit ball centred at $x_{i}$ (for the re-scaled metric), contradicting our assumption.

Thus we can apply Proposition 16 to construct  holomorphic maps $F_{i}:B_{i}\rightarrow \bC^{q}$ where $B_{i}\subset X_{i}$ is the set of points of (scaled) distance at most $\rho_{1}$ from $p$. Fix $h>0$ such that the $h$-neighbourhood of the set $\{\vert y\vert_{*}\leq 1/100\}$ in $\bC^{q}$ lies in the set $\{\vert y\vert _*\leq 1/2\}$. 
Pick a point $x\in D_{i}$ with $d_{i}(x,p)< \rho_{2}$ (which is possible for large $i$ by the definition of Gromov-Hausdorff convergence) and set $y=F_{i}(x)$ so $\vert y \vert_{*}\leq 1/100$. The fibre $F_{i}^{-1}(y)$ cannot meet the boundary of $B_{i}$ and so is a compact analytic subset. Since, as in the proof of Proposition 12, the line bundle $L^{k}$ is trivial over $B_{i}$ we see that the fibre is a finite set. It follows that the image  $F_{i}(D_{i})$ is an $(n-1)$-dimensional analytic set. More precisely the intersection with the ball of radius $h$ centred at $y$ is a closed $(n-1)$-dimensional analytic set and so its volume is bounded below, by the discussion in the proof of Proposition 14. The  bound on the derivative of $F_{i}$ gives a lower bound on the volume of $D_{i}\cap B_{i}$ and since $\rho_{1}<1$ this yields a lower bound on the volume of $D_{i}$ in the unit ball centred at $x_{i}$, for large $i$. This  gives our contradiction. 

\

{\bf Discussion}

The construction in Proposition 16 can be refined to show that  maps like $F_{i}$ can be chosen to  be embeddings. It is also possible to pass to a limit as $i\rightarrow \infty$, and thus get a local complex analytic model for the scaled limit space $Z'$. It is interesting to understand the algebro-geometric meaning of these. But here we have just constructed what we  need for our present purposes, and leave further developments for the future.

\subsection{Good tangent cones}

We consider a general scaled limit $Z'= \lim (X_{i}, a_{i} \omega_{i}) $ as above. We define a variant of the density function; for $z\in Z'$ and $\rho>0$
\begin{equation} V(i,z,\rho)=  \rho^{2-2n} {\rm Vol}\{ x \in D_{i}\subset X_{i}: d_{i}(x,z)\leq \rho\}, \end{equation} 
where the volume is computed using the metric $a_{i}\omega_{i}$.

 \begin{prop}
 There is a $c'>0$ such that for all $z\in {\cal D}(Z')$ and all $\rho>0$ we have
 $$  \liminf_{i} V(i,z, \rho)\geq c. $$
 \end{prop}
 This follows easily from Proposition 17.
 
 Now in the other direction we have
 \begin{prop} For any $z$ in $Z'$ which is not in $S_{2}(Z')$ there is a $\rho>0$ such that $V(i,z,\rho)$ is bounded.
\end{prop}
This follows immediately from Proposition  15.

\

The proof of Theorem 3 is a little complicated. For purposes of exposition we  will first prove a different but very similar  result.
\begin{prop}
Let $Z'$ be a limit space as above and let $K$ be a compact subset of $Z'$. Then
$S(Z')\cap K$ has capacity zero. 
\end{prop}

 The proof has seven  steps.

\

{\bf Step 1} By the Hausdorff dimension property we can find a countable collection $B_{\mu}$ of open balls of radius $r_{\mu}$ with $\sum r_{\mu}^{2n-3}$ arbitrarily small which cover $S_{2}\cap K$. We can suppose that the (1/10)-sized balls with the same centres are disjoint.  Let $U=\bigcup_{\mu} B_{\mu}$ so $U$ is open.

 {\bf Step 2} Let $J=K\cap (Z'\setminus U)$. Thus $J$ is a compact set and no point
of $J$ lies in $S_{2}$. Thus we can apply Proposition 19 at each point of $J$. Taking a finite covering we arrive at an open set $W$ such that $J\subset W$ and a $\rho>0$ such that the (2n-2)-volume of the intersection of $D_{i}$ with the $\rho$-neighbourhood of $\overline{W}$ is bounded by a fixed constant, for all $i$.

{\bf Step 3} The set $K\cap (Z'\setminus W)$ is compact and contained in $\bigcup_{\mu} B_{\mu}$. Thus we can find a finite sub cover, say by the $B_{\mu}$ for $\mu\leq M$. Then we can use the method of \cite{kn:DS} to construct a function $g_{1}$, equal to $1$ on $\bigcup_{\mu\leq M} B_{\mu}$, supported in a slightly larger set and with $L^{2}$ norm of the  derivative small. (Note that the point here is that the construction of \cite{kn:DS} will not work for infinite covers.)

{\bf Step 4} Now let $T$ be the intersection of the closure of ${\cal D}$ with $\overline{W}\cap K$. We want to show that $T$ has finite Minkowski measure. So we need to show that there  is an $M>0$ such that for any $\epsilon>0$ we can cover $T$ by at most $M \epsilon^{2n-2}$ balls of radius $\epsilon$. We may suppose that $\epsilon<\rho$ where $\rho $ is as in Step 2 above. By a standard argument we find a cover by $\epsilon$- balls such that the half-size balls are disjoint. Then the estimate on the number of balls follows from the upper and lower volume bounds for the $D_{i}$ (Propositions 18,19).

{\bf Step 5} By Proposition 1 we can  construct a cut-off $g_{2}$ with $L^{2}$ norm of the derivative small,  equal to $1$ on an open neighbourhood  $N$ of $T$ and supported in a somewhat larger (but arbitrarily small) neighbourhood.

{\bf Step 6} Now $g_{1}+ g_{2}$ is $\geq 1$ on the open set $N^{+}=N \cup \bigcup_{\mu\leq M} B_{\mu}$. Let $x$ be a point of $(S\setminus S_{2})\cap K$. Then if $x\in \overline{W}$ we have $x\in T$ so $x\in N$. If $x$ is not in $\overline{W}$ then $x$ is in $\bigcup_{\mu\leq M} B_{\mu}$. So, either way, $x$ is in $N^{+}$.
Write $\Sigma$ for the intersection of $S\cap K$ with $Z\setminus N^{+}$. Then $\Sigma$ is compact and, by the previous sentence, is contained in $S_{2}$.
Thus $\Sigma$ has Hausdorff codimension $>2$. As before, we can construct a cut-off $g_{3}$ for $\Sigma$ supported on a an arbitrarily small neighbourhod of $\Sigma$ and with derivative arbitrarily small in $L^{2}$. 
 
{\bf Step 7} Now $g_{1}+g_{2}+g_{3}$ is $\geq 1$ on all points of $S\cap K$. We choose a function $F(t)$ with $F(0)=0$ and with $F(t)=1$ if $t\geq 0.9$. Then set $g= F(g_{1}+g_{2}+g_{3})$.

\

Now we give the proof of Theorem 3, which is a variant of that above. Fix small $\delta>0$ and let $I$ be the interval
$[1-\delta, 1+\delta]$. Consider $Y\times I$ as embedded in the cone $C(Y)$ in the obvious way and identify $Y$ with $Y\times \{1\}$. To begin with we work on $Y$ so we cover $S(Y)$ by balls $B_{\mu}\subset Y$ as in Step 1 above and let $U=\bigcup_{\mu} B_{\mu}\subset Y$.
As in Step 2, we find an open set $W$ containing $Y\setminus U$ and a $\rho>0$ such that the volume of $D_{i}$ in the $\rho$-neighbourhood of $\overline{W}\times I$ is bounded. Then we find a finite set of balls $B_{\mu}$ which cover $ C(Y)\setminus W$. We construct a cut-off $g_{1}$ on $Y$ using this finite set of balls as in Step 3. Let $T$ be the intersection of the closure of $S(Y)\setminus S_{2}(Y)$ with $\overline{W}$. We want to show that the volume of the $\epsilon$-neighbourhood of $T$ is $O(\epsilon^{2})$. But in an obvious way this volume is controlled by the volume of the $\epsilon$-neighbourhood of $T\times I$ in $C(Y)$ for which we can argue as in Step 4 above. Then we construct a cut-off $g_{2}$ on $Y$ and the last steps  are just the same as above. 
\section{Proofs of Theorems 1 and 2}
\subsection{Global structure}

We want to define a notion of the complex structure on $Z$ being smooth in a neighbourhood of a point. We adopt a rather complicated definition, which is tailored to our needs. Given $p\in Z$, say the complex structure near $p$ is smooth if the following is true. There is a neighbourhood $\Omega_{\infty}$ of $p$ which is the limit of open sets $\Omega_{i}\subset X_{i}$ and  for some $k=k(p)$ there are sections $\sigma_{ i,0}, \dots , \sigma_{i,n-1}$ of $L^{k}\rightarrow X_{i}$ with the  properties below.
\begin{enumerate}
\item The $\sigma_{i,j}$ satisfy a fixed $L^{\infty}$ bound.
\item $\vert \sigma_{i,0}\vert\geq c>0$ on $\Omega_{i}$ for some fixed $c$.
\item The holomorphic maps $F_{i}:\Omega_{i}\rightarrow \bC^{n}$ given by
$\sigma_{i,j}/\sigma_{i,0}$ are homeomorphisms from $\Omega_{i}$ onto $B^{2n}\subset \bC^{n}$, satisfying a fixed Lipschitz bound.
\item The limit $F_{\infty}:\Omega_{\infty}\rightarrow B^{2n}$ is also a homeomorphism. 
\item There are holomorphic functions $P_{i}$ on $B^{2n}$ such that (in a smaller domain) the divisor $F_{i}(D_{i})$ is defined by the equation $P_{i}=0$ and the $P_{i}$ converge to a function $P_{\infty}$, not identically zero, as $i$ tends to $\infty$.
\end{enumerate}
Let $Z_{0}\subset Z$ be the set of points around which the complex structure is smooth. 
 By the nature of the definition this set is open in $Z$.
\begin{prop}
All points of ${\cal D}\cup R$ lie in $Z_{0}$.
\end{prop}
It follows that the complement $Z\setminus Z_{0}$ is a closed set of Hausdorff codimension at least $4$. This implies in particular that $Z_{0}$ is  connected.

The essential content of Proposition 21 is that ${\cal D}$ is contained in $Z_{0}$. The proof for $R$ is similar but easier. So suppose $p\in {\cal D}$ and  go back to the discussion in (2.5).  We have holomorphic equivalences
$$F_{i}: \Omega_{i}\rightarrow B^{2n}$$
where $\Omega_{i}\subset X_{i}$ lies between $B^{\sharp}(p,2)$ and $B^{\sharp}(p,1)$ and these maps satisfy a fixed Lipschitz bound. It follows that we can pass to the limit to get a Lipschitz map $F_{\infty}$ from a neighbourhood $\Omega_{\infty}$ of  $p$ in $Z$ to $B^{2n}$. The difficulty is that while we know the $F_{i}$ are injective we do not yet know the same for $F_{\infty}$. 

\begin{prop}
The map $F_{\infty}:\Omega_{\infty}\rightarrow B^{2n}$ is a homeomorphism.
\end{prop}
The essential thing is to show that the map is injective. So suppose we have
distinct $x,y\in \Omega_{\infty}$ with $F_{\infty}(x)=F_{\infty}(y)$.
Choose sequences $x_{i}, y_{i}$ in $X_{i}$ converging to $x,y$. Recall that we have locally non-vanishing sections $s$ of $L^{k}\rightarrow X_{i}$. Taking a subsequence, and possibly interchanging $x,y$ there is no loss in assuming that $\vert s(x_{i})\vert\geq \vert s(y_{i})\vert$ for all $i$. Applying Proposition 8 we can find some large $m$ such that there is for each large $i$ a section
$\tau$ of $L^{mk}$ such that $\vert \tau(x_{i})\vert =1$ and $\vert \tau(y_{i})\vert\leq 1/2$ (say).  Now set 
$$  f_{i}= \frac{\tau}{s^{m}}. $$
Then by construction $\vert f(x_{i}) \vert \geq 2 \vert f(y_{i})\vert $ and $\vert f(x_{i})\vert \geq c >0$ for some fixed $c$ (depending on $m$). Now consider the maps $(F_{i}, f_{i}):\Omega_{i}\rightarrow B^{2n}\times \bC$. On the one hand these satisfy a Lipschitz bound so we can pass to the limit and get a map $(F_{\infty}, f_{\infty})$ from $\Omega_{\infty}$ to
$B^{2n} \times \bC$ which separates $x,y$, by construction. On the other hand for finite $i$ the image of $(F_{i}, f_{i})$ is the graph of a holomorphic function $h_{i}$  on $B^{2n}$. The $h_{i}$ satisfy
a fixed $L^{\infty}$ bound so (taking a subsequence) we can suppose they converge to $h_{\infty}$. Clearly the image of $(F_{\infty}, g_{\infty})$ is the graph of $h_{\infty}$, which contradicts the fact that $(F_{\infty}, f_{\infty})$ separates $x,y$.

\

Proposition 22 now follows easily. For each $i$ the divisor is defined by a Weierstrass polynomial (19) and it is clear that (taking a subsequence) these have a nontrivial limit.  

\

{\bf Remark} One can give a different proof of Proposition 22 using the K\"ahler-Einstein equations, similar to Proposition 2.5  in \cite{kn:CDS1}. We will give a version of this argument in the sequel. The advantage of the proof given is that it does not use the K\"ahler-Einstein equations and potentially extends to other situations.

\

Given this background we can proceed, following the arguments in \cite{kn:DS}, to analyse the global structure of $Z$.
\begin{itemize}
\item The arguments in \cite{kn:DS} go over without any change to show that there is some large $m$ with the following effect. Use the $L^{2}$ norm to define a metric on $H^{0}(X_{i}, L^{m})$ and pick an orthonormal basis. Then we get projective embeddings   $T_{i}:X_{i}\rightarrow \bC\bP^{N}$, canonically defined up to unitary transformations. The images converge to a normal
variety $W\subset \bC\bP^{N}$ and  there is a Lipschitz homeomorphism $h:Z\rightarrow W$ compatible with the Gromov-Hausdorff convergence, in the sense made precise in \cite{kn:DS}. 
\item Likewise the set $Z_{0}\subset Z$ maps under $h$ to the smooth part of $W$. This is clear from the definition of $Z_{0}$. (More precisely,
that definition shows that a neighbourhood of $p\in Z_{0}$ is embedded smoothly in projective space by sections of $L^{k(p)}$. But once we know that the images of $Z$ stabilise for large $k$ we can take a fixed power $m$.)
\item It is also clear that $h(Z_{0})$ is exactly equal to the smooth part of $W$.
\item Taking a subsequence, we can suppose that the $T_{i}(D_{i})$ converge to an $n-1$ dimensional algebraic cycle $\Delta= \sum_{a} \nu_{a} \Delta_{a}$ in $\bC\bP^{N}$. The fact that $W$ is normal implies that the intersection of each $\Delta_{a}$ with the singular set of $W$ is a proper algebraic subset of $\Delta_{a}$.
Then $\Delta$ is a Weil divisor in $W$.
\item Since the volume of  the $D_{i}$ is a fixed number  the lower  bound on densities (Proposition 17) implies a fixed bound on the Minkowski measure $m(D_{i})$. Then, as in the discussion at the end of (2.1), we get (after perhaps passing to a subsequence) a limit set $D_{\infty} \subset Z$ which is a closed set with $m(D_{\infty})\leq M$. Thus $D_{\infty}$ has finite codimension 2 Hausdorff measure. 
\item We have ${\cal D}\subset D_{\infty} \subset S(Z)$. This is clear from the discussion in (2.2). 
\item The homeomorphism $h:Z\rightarrow W$ maps $D_{\infty}$ onto ${\rm supp} \Delta$. This follows from the compatability between the algebraic and Gromov-Hausdorff convergence. 
\item The complement $S_{2}= S\setminus {\cal D}$ has Hausdorff dimension at most $2n-4$. Thus the same is true for $D_{\infty}\setminus {\cal D}$. The Lipschitz property of $h$ implies that $ {\rm supp}\ \Delta  \setminus
h({\cal D})$ has Hausdorff dimension at most $2n-4$. In particular
$h({\cal D})$ is dense in ${\rm supp} \Delta $.

\end{itemize} 

All of the above discussion would apply equally well (with suitable hypotheses) in a more general situation where $L$ is not  equal to the anti-canonical bundle. We will now develop results specific to the \lq\lq Fano case'', which is our concern in this paper.  For simplicity we assume that $\lambda=1$. (The general case would be exactly the same, since we can take $k$ to run through multiples of the fixed number $\lambda$.)

For each $i$, the divisor $D_{i}\subset X_{i}$ is defined by a section
$S_{i}$ of $K^{-1}$. We have a hermitian metric $h_{i}$ on $L=K^{-1}$. The K\"ahler-Einstein equation takes the form
\begin{equation}  \omega_{i}^{n}= (S_{i}\wedge \overline{S_{i}})^{-1} \vert S_{i} \vert^{2\beta_{i}}, \end{equation}
and we can normalise so that the $L^{2}$ norm of $S_{i}$ (defined using the metric $h_{i}$) is $1$. 
 Then from (7) we get a fixed upper bound on the $L^{\infty}$ norm of $S_{i}$. 

Now consider the situation around a point of $Z_{0}$ in the charts $\Omega_{i}$ occurring in the definition.  We have local co-ordinates $z_{a}$ and a K\"ahler potential $\phi_{i}= k^{-1} \log \vert \sigma_{0,i}\vert^{2}$ which is bounded above and below. Let $\Theta=dz_{1}\dots dz_{n}$. We can write $\sigma_{0,i}= V_{i}^{k} \Theta^{-k}, S_{i}= U_{i} P_{i} \Theta^{-1}$ where $U_{i},V_{i}$ are non-vanishing holomorphic functions on the ball, and $P_{i}$ is the given local defining function for $D_{i}$.  Then the equation becomes (dropping the index $i$)
\begin{equation}
\det(\phi_{a\overline{b}}) = \vert U P \vert^{2\beta-2} \vert V\vert^{-2\beta} e^{-\beta \phi}. \end{equation}
\begin{prop}
In such co-ordinates there is $C$, independent of $i$, such that on a fixed interior ball
$C^{-1}< \vert U_{i}\vert, \vert V_{i}\vert < C$.
\end{prop}

 Consider $\mu>0$, an open set $G$ in $B^{2n}$  and the integral
$$  I_{\mu, G}= \int_{G} \vert S_{i}\vert^{2\mu} \omega_{i}^{n}. $$
The $L^{\infty}$ bound on $S_{i}$ gives an upper bound on $I_{\mu, G}$, independent of $i$. On the other hand the normalisation $\Vert S_{i}\Vert_{L^{2}}=1$ and an analytic continuation argument show that, for any $\mu,G$, the integral $I_{\mu, G}$  has a strictly positive lower bound, independent of $i$. Now, using the K\"ahler-Einstein equation, we can write
$$  I_{\mu, G}= \int_{G} \vert P U\vert^{2\beta+2\mu-2} \vert V\vert^{-2\beta-2\mu} e^{-(\mu+ \beta) \phi} \Theta \wedge \overline{\Theta} , $$
where again we temporarily drop the index $i$. 
First take $\mu=1-\beta$. Then the integral becomes
  $$  I= \int \vert V\vert^{-2} e^{-\phi} \Theta \wedge \overline{\Theta}, $$
 Since $\phi$ is bounded above and below this is equivalent to the integral of $\vert V\vert^{-2}$, with respect to the Euclidean measure. Begin by taking $G$ to be the whole unit ball. Thus we have an $L^{2}$ bound on the holomorphic function $f=V^{-1}$ over the ball which gives an $L^{\infty}$ bound $\vert  f\vert \leq C$ over an interior ball. We claim that on the other hand we must have some lower bound $\vert f \vert \geq C^{-1}$ over an interior ball, say $(1/4) B^{2n}$. To see this we argue by contradiction: if not there is a violating sequence $f_{j}$ with points $z_{j}\in (1/4)B^{2n}$ such that $\vert f_{j}(z_{j})\vert \rightarrow 0$. Taking a subsequence, we can suppose that the $z_{j}$ have limit $z$ and the $f_{j}$ converge to a limit $f_{\infty}$ in $C^{\infty}$ on compact subset of $B^{2n}$. We must have $f_{\infty}(z_{\infty})=0$.
If $f_{\infty}$ is not identically zero then, by basic complex analysis, we get a nearby zero of $f_{j}$ for large $j$ in contradiction to our hypothesis. On the other we know that the integral of $\vert f_{j}\vert^{2}$ over a small ball centred at $z_{\infty}$ has a strictly positive lower bound, so the limit cannot be identically zero. 

This argument gives the required upper and lower bounds on $\vert V\vert$. Now take $\mu=2-\beta$. Then $I_{\mu, G}$ is equivalent to the integral over $G$ of $\vert U \vert^{2} \vert P\vert^{2}$. Re-instate the index $i$ so $P=P_{i}$ which have non-trivial limit $P_{\infty}$. Away from the zero set of $P_{\infty}$ we can apply the argument above to get upper and lower bounds on $\vert U_{i}\vert$, for large $i$. Then we can extend these over the zero set by applying the Cauchy integral formula to the restrictions of $U_{i}, U_{i}^{-1}$ to suitable small discs.

\

\

Given Proposition 23, we can pass to the limit  $i=\infty$ (working over a slightly smaller ball). The upper and lower bounds on $\vert U_{i}\vert, \vert V_{i}\vert$ mean that we can take limits and get holomorphic functions $U_{\infty}, V_{\infty}$ with the same bounds. Then $U_{\infty}P_{\infty} \Theta$ is a local section of $K_{W}^{-1}$ defining the divisor $\Delta$.
While this is a local discussion the local sections obviously glue together to define a global section $S_{\infty}$ of $K_{W}^{-1}$ which defines $\Delta$ (in the sense that $\Delta$ is the closure of the zero set of $S_{\infty}$ over the smooth part of $W$, as discussed in Section 1).  

\
 Now to see that $(W,(1-\beta_{\infty})\Delta)$ is KLT we write 
 $$  (S\wedge \overline{S})^{-1}= (\sigma \wedge \overline{\sigma})^{-1/m} \vert \sigma \vert^{2/m} \vert S \vert^{-2}. $$
(Here we drop the subscript $\infty$.)  So the K\"ahler-Einstein equation gives
$$  (\sigma \wedge \overline{\sigma})^{-1} \vert S\vert^{2\beta-2}=  \omega^{n} \vert \sigma \vert^{-2/m}. $$
The KLT condition is that the integral over the ball  of $(\sigma \wedge \overline{\sigma})^{-1} \vert S\vert^{2\beta-2}$  is finite but by the identity above this is the integral of $\omega^{n} \vert \sigma \vert^{-2/m}$ which is clearly finite (since $\sigma$ has no zeros).
and the finiteness of the integral is clear (since $\vert \sigma\vert$ is bounded below).  It also follows immediately from the set-up that the limiting metric is a weak conical K\"ahler-Einstein metric, as defined in Section 1.  This completes the proof of Theorem 1.

\subsection{Smooth limits: Proof of Theorem 2}

We now turn to Theorem 2. Of course what we prove is a local result: if $p$ is  a smooth point of $W$ and of the divisor $\Delta$ then, near $p$ the limiting metric is a metric with cone singularities along $\Delta$ and cone angle $2\pi\beta_{\infty}$. Recall from the introduction that, for our purposes, this means in the sense of the class of metrics defined in \cite{kn:Dona11}, so we briefly review this definition. Work in local co-ordinates $(u,v_{1},\dots, v_{n-1})$ and let $\omega_{(\beta)}$ be the standard metric with a cone singularity of angle $2\pi \beta$ along $u=0$. Let $r=\vert u\vert^{\beta}$ and $\theta$ be the argument of $u$ (in the obvious sense). Then we have an orthonormal basis for the (1,1) forms (in the model metric $\omega_{(\beta)}$)  
$$  \epsilon \wedge \overline{\epsilon},\  \epsilon \wedge d\overline{v}_{i},\  \overline{\epsilon} \wedge dv_{i}, \ dv_{i}\wedge d\overline{v}_{j}, $$
where $\epsilon= dr+ i\beta^{-1} rd\theta$. These are well-defined forms on $\{u\neq 0\}$, even though the co-ordinate $\theta$ is only locally defined. We say that a $(1,1)$ form  is in ${\cal C}^{,\alpha,\beta}$ if its components in this frame are $C^{,\alpha}$ functions with respect to the metric $\omega_{(\beta)}$. Likewise for a function in ${\cal C}^{,\alpha,\beta}$. We say that a function $\phi$ is in ${\cal C}^{2,\alpha,\beta}$ if $\phi$ and $\dbd \phi$ are in ${\cal C}^{,\alpha,\beta}$. We say that a metric $\omega$ has cone singularities with cone angle $2\pi \beta$ along $\{u=0\}$ if 
\begin{itemize}
\item $\omega$ is uniformly equivalent to $\omega_{(\beta)}$, so
$$  C^{-1} \omega_{(\beta)}\leq \omega\leq C \omega_{(\beta)}$$
for some $C$;
\item $\omega=\omega_{(\beta)}+ \dbd \phi$ where $\phi\in {\cal C}^{2,\alpha,\beta}$.
\end{itemize}

\

  Now  go back to the converging sequence of metrics $\omega_{i}$ on $X_{i}$, working in local co-ordinates $(u,v_{1}, \dots, v_{n-1})$. There is no loss of generality in supposing that, near $p$ each $D_{i}$ is given by $u=0$ (since we can always adjust the co-ordinates to achieve this, on a small ball). For any angle $\beta$ let $\omega_{(\beta)}$ be the standard cone metric on $\bC^{n}= \bC_{\beta}\times \bC^{n-1}$, as above.  We first establish a uniform bound.
\begin{prop}
There is a constant $C$ such that $$C^{-1} \omega_{(\beta_{i})} \leq \omega_{i}\leq C \omega_{(\beta_{i})}. $$
\end{prop}

To prove Proposition 24, drop the index $i$ and write $\beta_{i}=\beta$ and  $\tomega_{i}=\omega$. We have $\omega=\dbd \phi$ where $0\leq \phi\leq C_{0}$. In the set up of (3.1) the K\"ahler-Einstein equation is
$$   \omega^{n}=  \vert u\vert^{2\beta-2} \vert Q\vert^{2} e^{-\beta \phi} \Theta \wedge \overline{\Theta}, $$
where $Q$ is a holomorphic function on the ball with $C_{1}^{-1}\leq \vert Q\vert\leq C_{1}$. (That is,  $Q= U^{\beta-1} V^{-\beta}$, in the notation of (3.1)---we can take fractional powers since the ball is simply connected.) Here $\Theta=du dv_{1}\dots dv_{n-1}$. However changing $\phi$ by subtracting
$\beta^{-1}\log \vert Q\vert^{2}$ we can suppose that in fact $Q=1$ so our equation is

\begin{equation}\omega^{n}=    e^{-\beta \phi} \omega_{(\beta)}^{n}, \end{equation}

since $\omega_{(\beta)}^{n}= \vert u\vert^{2\beta-2} \Theta \wedge \overline{\Theta}$.
We know that $\omega\geq C_{3}^{-1} \omega_{{\rm Euc}}$, where $\omega_{{\rm Euc}}$ is the standard Euclidean metric on $\bC^{n}$ (by item (3) in the definition at the beginning of (3.1)). To establish Proposition 24 it suffices to show that
\begin{equation}   \omega \geq C_{4}^{-1} \omega_{(\beta)}, \end{equation}
 for by (26) the ratio of the volume forms $\omega^{n}, \omega_{(\beta)}^{n}$ is bounded above and below and the lower bound (27) gives an upper bound.

To prove (27), let $h$ be the square root of the trace of $\omega_{(\beta)}$ with respect to $\omega$, so
\begin{equation} h^{2} \omega^{n}= \omega_{(\beta)} \wedge \omega^{n-1}.\end{equation}

We claim that on a smaller ball, say $(1/2) B^{2n}$, we have an integral bound
\begin{equation}  \int_{(1/2) B^{2n}}  h^{2} \omega^{n} \leq C_{5}. \end{equation}
This follows from the Chern-Levine-Nirenberg argument. Let $\chi$ be a standard function of compact support on $B^{2n}$ equal to $1$ on $(3/4) B^{2n}$, say. Then we have $\dbd \chi \leq C_{6} \omega_{{\rm Euc}} \leq C_{3} C_{6} \omega_{(\beta)}$.
Then $$\int_{(3/4)B^{2n}} \omega_{\beta}\wedge \omega^{n-1}\leq \int_{B^{2n}} \chi \omega_{\beta}\wedge \omega^{n-1}. $$
Integrating by parts, the right hand side is
$$  \int_{B^{2n}}\phi ( \dbd \chi) \omega_{\beta} \omega^{n-2}, $$
which is bounded by $$C_{0} C_{6} \int_{B^{2n}} \omega_{\beta}^{2} \omega^{n-2}.$$
Just as in Proposition 15 above, we continue with a nested sequence of cut-off functions to interchange powers of $\omega_{\beta}, \omega$ and arrive at a bound (29). 
The function $h^{2}$ can also be regarded as $\vert \nabla f\vert^{2}$ where
$f$ is the identify map from the ball with metric $\omega$ to the ball with metric $\omega_{(\beta)}$. The fact that $\omega$ has positive Ricci curvature implies that away from the singular set we have $\Delta h \geq 0$. In fact we have \begin{equation} \Delta h^{2} =\Delta (\vert \nabla f\vert^{2}) \geq 2 \vert \nabla \nabla f \vert^{2} \end{equation}
where $\nabla f \in \Gamma(T^{*}\otimes T)$ is the identity endomorphism of the tangent bundle and the covariant derivative is that formed using $\omega$ on one factor and $\omega_{(\beta)}$ on the other. This can be seen as an instance of the Chern-Lu inequality, or more simply as the Bochner formula, since away from the singular set we can take local Euclidean co-ordinates for $\omega_{(\beta)}$. Using the fact that $\vert \nabla \nabla f\vert\geq \vert \nabla \vert \nabla f \vert \vert$ we get $\Delta h\geq 0$. 
\begin{lem} The inequality $\Delta h \geq 0$ holds in a weak sense across the singular set.  That is, for any compactly supported smooth non-negative test function $F$ we have
\begin{equation}   \int (\nabla h, \nabla F) \leq 0. \end{equation}
\end{lem}

It is here that we use the fact that the cone angles match up. This implies, by the definition of a cone singularity, that $h$ is a bounded function on the complement of the divisor $\{u=0\}$, say $h\leq M$. Thus any given metric $\omega$ is bounded above and below by {\it some} multiples of $\omega_{(\beta)}$.  Let $\sigma$ be a non-negative function of compact support on the intersection of the ball with $\{u\neq 0\}$. Then $$ \nabla(\sigma^{2} h). \nabla h = \vert \nabla(\sigma h)\vert^{2} - h^{2} \vert \nabla \sigma\vert^{2}. $$
Since $$\int \nabla(\sigma^{2} h ). \nabla h = -\int(\sigma^{2} h \Delta h) \leq 0, $$
we have
$$  \int \vert \nabla(\sigma h)\vert^{2} \leq \int h^{2} \vert \nabla \sigma \vert^{2}\leq M^{2} \int \vert \nabla \sigma \vert^{2}. $$
By a construction like that in (2.1) we can choose functions $\sigma_{j}$ equal to $1$ in intersection of the half-ball $\frac{1}{2} B^{2n}$ with $\{\vert u\vert\geq \epsilon_{j}\}$ where $\epsilon_{j}\rightarrow 0$ but with the $L^{2}$ norm of $ \nabla \sigma_{j} $ bounded by a fixed constant. (It suffices to do this using the model metric $\omega_{(\beta)}$ since, as we have noted above, the given metric $\omega$ is bounded above and below by {\it some} multiples of $\omega_{(\beta)}$.) Taking the limit as $j\rightarrow \infty$ we see that
$$  \int_{\frac{1}{2}B^{2n}} \vert \nabla h \vert^{2} <\infty. $$
By obvious arguments, it suffices for our purposes to prove (31) for functions $F$ supported in $\frac{1}{2} B^{2n}$. Thus we can suppose that $F \nabla \sigma_{j}$ is supported in $\{ u\leq \epsilon_{j}\}$ and that the $L^{2}$ norm of $F \nabla \sigma_{j}$ tends to $0$ as $j\rightarrow \infty$. Now
$$  \int \sigma_{j} \nabla F. \nabla h = -\int F \nabla \sigma_{j}.\nabla h - \int F \sigma_{j} \Delta h \leq \Vert F \nabla \sigma_{j}\Vert_{L^{2}} \Vert \nabla h\Vert_{L^{2}}. $$
Taking the limit as $j\rightarrow \infty$ we see that 
$$  \int \nabla F.\nabla h \leq 0. $$

With this Lemma in place we return to the proof of Proposition 24. We use the Moser iteration technique, starting with the bound on the integral of $h^{2}$ and the inequality $\Delta h \geq 0$  to obtain an $L^{\infty}$ bound on $h$.  Here of course we use the fact that there is a uniform Sobolev inequality for $\omega$. The Euclidean ball  $(1/2)B^{2n}$ in our co-ordinates contains the metric ball (in the metric $\omega$ of radius $(2 C_{3})^{-1}$ centred at $p$ (since $\omega\geq C_{3}^{-1}\omega_{{\rm Euc}}$). Then Moser iteration gives a bound on $h$ over the metric ball of radius $(4 C_{3})^{-1}$ say. Finally this metric ball contains a Euclidean ball of a definite size. This can be seen from the argument of \cite{kn:CDS1} Proposition 2.4. (In fact that argument shows that the identity map from $(B^{2n}, \omega_{{\rm Euc}})$ to $(B^{2n}, \omega)$ satisfies a fixed H\"older bound with exponent $\beta$.) This completes the proof of Proposition 24.

To deduce Theorem 2 from Proposition 24 we need some way of improving the estimates. In this general area, one approach is to try to estimate higher derivatives. For  metrics with cone singularities, this was achieved by Brendle \cite{kn:Brendle} assuming that $\beta<1/2$. Another approach in this general area is furnished by the Evans-Krylov theory and an analogue of that theory for metrics with cone singularities was developed by Calamai and Zheng \cite{kn:CZ}, assuming that $\beta<2/3$. A general result of Evans-Krylov type is stated in [13] where two independent proofs are given, but at the time of writing we have had difficulty in following the one of these proofs that we have, so far, studied in detail. Partly for this reason, and partly because it has its own interest,  we give a different approach (to achieve what we need) below. Our approach is related to ideas of Anderson \cite{kn:A} in the standard theory, and we have other results, which apply under slightly different hypotheses, using a line of argument closer to Anderson's, and which we will give elsewhere.

\

Write $\omega=\dbd \phi$ for the limit of $\omega_{i}$ as $i$ tends to infinity. This is a smooth K\"ahler-Einstein metric outside $\{u=0\}$ and is a local representation of  the metric on the Gromov-Hausdorff limit. Proposition 24 and the $C^{1}$ estimate (8) imply that $\phi$ satisfies a Lipschitz bound with respect to the standard cone metric $\omega_{(\beta)}$. 
 
 For $\epsilon>0$ let $T_{\epsilon}: \bC^{n}\rightarrow \bC^{n}$ be the linear map $T_{\epsilon}(u,\underline{v})= \epsilon^{1/\beta} u, \epsilon^{-1} \underline{v})$.
Suppose $\epsilon_{j}$ is any sequence with $\epsilon_{j}\rightarrow 0$ and define a sequence $\omega^{(j)}$ of re-scalings by 
\begin{equation}\omega^{(j)}= \epsilon^{-2} T_{\epsilon}^{*} (\omega).\end{equation}
The scaling behaviour of the model $\omega_{(\beta)}$ implies that we still have uniform upper and lower bounds $C^{-1} \omega_{(\beta)}\leq \omega^{(j)}\leq C \omega_{(\beta)}$ so (perhaps taking a subsequence) we get $C^{\infty}$ convergence on compact subsets of $\{u\neq 0\}$ to a limit $\omega^{(\infty)}$. The K\"ahler-Einstein equation satisfied by $\phi$ and the fact that $\phi$ satisfies a Lipschitz bound with respect to the $\omega_{(\beta)}$ metric implies that $(\omega^{(\infty)})^{n}= \kappa \omega_{(\beta)}^{n}$ for some constant $\kappa>0$ which without loss of generality we can suppose is $1$. On the other hand the uniform bound of Proposition 24 implies that this limit $\omega^{(\infty)}$ represents the metric on a tangent cone to $(Z,\omega)$ at the given point $p$. 

\begin{prop} Suppose $\omega'$ is a K\"ahler metric defined on the complement of the divisor $\{u=0\}$ in $\bC^{n}$ such that
\begin{itemize}
\item The volume form of $\omega'$ is the same as that of $\omega_{(\beta)}$;
\item there are uniform bounds $$  C^{-1} \omega_{(\beta)}\leq \omega'\leq C \omega_{(\beta)};$$
\item $(\bC^{n},\omega')$ is a metric cone with vertex $0$.
\end{itemize}
Then there is a complex linear isomorphism $g:\bC^{n}\rightarrow \bC^{n}$ preserving the subspace $\{u=0\}$ such that $\omega'=g^{*}(\omega_{(\beta)})$.
\end{prop}

In our context this implies that any tangent cone at a point of $\Delta$ is $\bC_{\beta}\times \bC^{n-1}$, in other words all points of $\Delta$ lie in ${\cal D}_{1}$. Conversely it is clear from the Weierstrasse representation (19) that at points of ${\cal D}_{1}$ the limiting divisor is smooth so we are in the situation considered here.

The hypothesis that $(\bC^{n},\omega')$ is a cone means that there is a radial vector field $X$, initially defined over $\{u\neq 0\}$. In the model case of $\omega_{(\beta)}$ the corresponding vector field is
$$  X_{(\beta)}= {\rm Re}\left( \beta^{-1} u \partial_{u} + \sum_{i=1}^{n-1} v_{i} \partial_{v_{i}} \right). $$
We want to show that, after perhaps applying a linear transformation we have $X=X_{(\beta)}$. Write 
$$ X={\rm Re} \left( a_{0} \partial_{u}+ \sum_{i=1}^{n-1} a_{i} \partial_{v_{i}}\right) , $$
for functions $a_{i}$. The K\"ahler condition means that the $a_{i}$ are holomorphic.  
The length of $X$, computed in the metric $\omega'$ is equal to the distance to the origin. By the uniform upper and lower bounds we get
\begin{equation}  \vert u\vert^{2\beta-2} \vert a_{0}\vert^{2} +\sum_{i=1}^{n-1} \vert a_{i}\vert^{2} \leq C'( \vert u \vert^{2\beta} + \sum_{i=1}^{n-1} \vert v_{i}\vert^{2}), \end{equation}
for another constant $C'$. This implies that for each $i\geq 0$ we have $$\vert a_{i}\vert \leq C'' \vert (u,\underline{v})\vert^{1+\beta} .$$
First, by  Riemann's removal singularities theorem this means  that the $a_{i}$ extend holomorphically over $u=0$. Second, since $\beta<1$ this means that the $a_{i}$ are linear functions. Also we see  (taking $u=0$ in (33)) that $a_{0}$ vanishes when $u=0$. Thus the vector field is defined by an endomorphism $A:\bC^{n}\rightarrow \bC^{n}$, preserving the subspace $\{u=0\}$. Write $A_{(\beta)}$ for the endomorphism corresponding to $X_{(\beta)}$.

For $\tau\in \bC$ consider the vector field ${\rm Re}(\tau X)$ and the one parameter group $f_{\tau}$ of holomorphic transformations so generated. Thus $f_{\tau}$ takes a point of distance $1$ from the origin to a point of distance
  $\exp({\rm Re}(\tau))$ from the origin. Using the equivalence of metrics again, this translates into the statement that there is a constant $C''$ such that for all $z\in \bC^{n}$ with $\vert z \vert =1$ we have
\begin{equation} (C'')^{-1}\exp({\rm Re}(\tau)) \leq [ \exp(\tau A) z ] \leq C'' \exp({\rm Re}(\tau))\end{equation}
where $[(u,\underline{v})]= \left( u^{2\beta} + \vert \underline{v}\vert^{2}\right)^{1/2}$.
Let $A_{0}$ be the restriction of $A$ to $\{u=0\}$ and apply this  to an eigenvector $v$ of $A_{0}$. We see that the eigenvalue must be $1$. Suppose that $A_{0}$ is not the identity. Then we can write $A_{0}= 1+ N$ where
for some $k\geq 1$ we have $N^{k}\neq 0$ but $N^{k+1}=0$. Then
$$  \exp(\tau A_{0}) = \exp(\tau) ( \tau^{k} N^{k}/k!+ \dots), $$
and if we take $v$ so that $N^{k} v\neq 0$ we see that this violates (34) for large $\tau$. So $A_{0}$ is the identity. Arguing in the same way we see that the there must be an eigenvalue $\beta^{-1}$ of $A$. So $A$ is conjugate to $A_{(\beta)}$ by a map preserving $\{u=0\}$ and without loss of generality we can suppose that $X=X_{(\beta)}$.

To complete the proof of Proposition 25, we consider the function $h$ defined as the square root of the trace of $\omega_{(\beta)}$ with respect to $\omega'$, just as in the proof of Proposition 24. This is defined on the complement of the singular set $\{u=0\}$ and is invariant under the real 1-parameter subgroup generated by $X=X_{(\beta)}$. As before  we have $\Delta h\geq 0$. Let ${\cal A}$ be the \lq\lq annulus'' consisting of points $z\in \bC^{n}$ with distance (in the metric $\omega'$) between $1$ and $2$ from the origin. Let $\chi$ be a non-negative  function on a neighbourhood of ${\cal A}$, equal to $1$ outside a small neighbourhood of $\{u=0\}$ and vanishing near to $\{u=0\}$. Then, much as in the proof of Lemma 2, 
$$ - \int_{{\cal A}} \chi^{2} h \Delta h = \int_{{\cal A}} \vert \nabla (\chi h )\vert^{2} - h^{2} \vert \nabla \chi\vert^{2} - \int_{\partial{\cal A}} \chi^{2} h (\nabla h.n), $$
where $n$ is the unit normal vector to the boundary.
The fact that $h$ is invariant under the flow of $X$ implies that $(\nabla h.n)=0$ so the boundary term vanishes and taking the limit over a sequence of functions $\chi$ we see that the integral of $\vert \nabla h\vert^{2}$ vanishes, so $h$ is a constant. Going back to (30), we see that $\omega$ is covariant constant with respect to $\omega_{(\beta)}$. This easily implies that $\omega$ is isometric to $\omega_{(\beta)}$ by a complex linear transformation.

An alternative  approach to this last part of the proof is to apply the maximum  principle to $h$. The only difficulty is that the maximum might be attained at a singular point $q=(u,0)$ with $u\neq 0$. We can work by induction on $n$ so we suppose we have proved Proposition 25 in lower dimensions. If we blow up
at this point $q$, we arrive at a K\"ahler Ricci flat metric cone $(C(Y)\times \bC^{n-k},\omega_{\phi}) \;$ which
is quasi isometric to $(\bC_{\beta}\times \bC^{n-1},\omega_{\beta}).\;$ By Cheeger-Colding theorem,
we know that $n-k \geq 1.\;$  The splitting of $C^{n-k}$ is represented by $n-k$ dimensional complex holomorphic
vector field with constant length w.r.t. $\omega_{\phi}.\;$ By quasi isometry, this set of holomorphic vector fields
is bounded in $(\bC_{\beta}\times \bC^{n-1},\omega_{\beta}).\;$ It follows these are constant holomorphic
vector fields in $\bC^{n}.\;$ It follows that $(C(Y)\times \bC^{n-k},\omega_{\phi}) \;$ splits holomorphic and
isometrically so that $(C(Y), \omega_{\phi}\mid_{C(Y)})$ is Ricci flat, quasi-isometry to $(\bC_{\beta}\times \bC^{k-1}, \omega_{\beta})$ in  $\bC_{\beta}\times \bC^{k-1}.\;$ By induction, we know that $(C(Y), \omega_{\phi}\mid_{C(Y)})$
is equivalent to $(\bC_{\beta}\times \bC^{k-1}, \omega_{\beta}).\;$ It follows that $(C(Y)\times \bC^{n-k},\omega_{\phi}) \;$
is standard.  Then (using also the arguments below) we can prove that $\omega'$ is $C^{,\alpha,\beta}$ near $q.\;$ Once we know this it is straightforward to see that the maximum principle can be applied.

We now go back to the metrics $\omega_{i}$, before taking the limit. For $\epsilon\in (0,1)$ let $$\omega_{i,\epsilon}= \epsilon^{-2} T^{*}_{\epsilon}(\omega_{i}). $$
We consider this as a metric defined on the unit ball $B^{2n}$. 
\begin{prop}
Given any $\zeta>0$ we can find $\epsilon(\zeta)$ such that for $\epsilon\leq \epsilon (\zeta)$ we can find a linear map $g_{\epsilon}$ preserving the subspace $\{u=0\}$ such that for  $i\geq i(\epsilon)$ we have 
\begin{equation} (1-\zeta) g_{\epsilon}^{*}(\omega_{(\beta_{i})})\leq \omega_{i,\epsilon}\leq (1+\zeta) g_{\epsilon}^{*}(\omega_{(\beta_{i})}), \end{equation}
throughout $\frac{1}{2} B^{2n}$.
\end{prop}
The argument is similar to that in the proof of Proposition 24. Let $h_{i,\epsilon}$ be the square root of the trace of $\omega_{(\beta_{i})}$ with respect to $\omega_{i,\epsilon}$. The uniform Lipschitz bound on the K\"ahler potential and the K\"ahler-Einstein equation (26) imply that the ratio of the volume forms $\omega_{i,\epsilon}^{n}, \omega_{(\beta_{i})}^{n}$ is within $c\epsilon$ of $1$, for a fixed constant $c$. Thus if $\lambda_{a}$ are the eigenvalues of $\omega_{(\beta_{i})}$ with respect to $\omega_{i,\epsilon}$ we have
$$   h^{2}_{i,\epsilon}= \sum \lambda_{a}, \prod \lambda_{a}\geq (1-c\epsilon). $$
Thus $$  h^{2}_{i,\epsilon}\geq n(1-c\epsilon)^{1/n}\geq n-c'\epsilon, $$
say, for a suitable fixed $c'$. Set $f_{i,\epsilon}= h^{2}_{i,\epsilon}-n+c'\epsilon$, so $f_{i,\epsilon}\geq 0$ and $\Delta f_{i,\epsilon}\geq 0$. Suppose we know that, for some $\theta$,
the $L^{2}$ norm of $f_{i,\epsilon}$ over $B^{2n}$ is bounded by $\theta$. Then Moser iteration implies that $f_{i,\epsilon}\leq K \theta$ over the half-sized ball, for a fixed computable $K$. Thus $\sum \lambda_{a}\leq n+ K\theta$ and $\prod \lambda_{a} \geq 1-c'\epsilon$. It follows by elementary arguments that $1-\zeta\leq \lambda_{a}\leq 1+\zeta$ provided that $\theta, \epsilon$ are sufficiently small compared with $\zeta$. 

The problem then is to show that for any $\theta$ we can make the $L^{2}$ norm of $f_{i,\epsilon}$ less than $\theta$, by taking $\epsilon$ small and then $i\geq i(\epsilon)$. Of course to do this we may need to apply a linear transformation $g_{\epsilon}$, as in the statement of Proposition 26. Clearly it is the same to get the $L^{2}$ bound on  $h^{2}_{i,\epsilon}- n$. First we work with the limiting metric $\omega_{\infty}$, writing $h_{\infty,\epsilon}$ in the obvious way. Since $h_{\infty,\epsilon}$ satisfies a fixed $L^{\infty}$ bound, it suffices to show that as  $\epsilon$ tends to zero the functions $h^{2}_{\infty, \epsilon}-n$ converge to zero uniformly on compact subsets of $\{u\neq 0\}$. But this follows from Proposition 25 (after applying linear transformations). Now the convergence of $\omega_{i}$ to $\omega_{\infty}$, again uniformly on compact subsets of $\{u\neq 0\}$ gives our result.\\


We now move on to  complete the proof of Theorem 2. Given any fixed $\zeta>0$ we can without loss of generality suppose (by the result above) that over the unit ball $B^{2n}$ our metrics $\omega_{i}$ satisfy \begin{equation}(1-\zeta)\omega_{(\beta_{i})}\leq \omega_{i}\leq (1+\zeta) \omega_{(\beta_{i})}. \end{equation}

 For $q\in \bC^{2n}$ and $\rho>0$ write
$B^{\beta_{i}}(q,\rho)$ for the ball of radius $\rho$ and centre $q$, in the metric $\omega_{(\beta_{i})}$.
\begin{lem}
There are $C, \zeta_{0}>0$ and $ R>2$ such that, if $\omega_{i}$ satisfies (36) on $ B^{2n}$ with $\zeta=\zeta_{0}$, then if $q$ is any point of $B^{2n}$ and $\rho$ is such that
$B^{\beta_{i}}(q, R \rho)\subset B^{2n}$ then on $B^{\beta_{i}}(q,2\rho)$ we can write $\omega_{i}=\omega_{(\beta_{i})}+\dbd \psi$ where
$[\psi]_{\alpha}\leq C \rho^{2-\alpha}$.
\end{lem}
Here $[\ ]_{\alpha}$ denotes the H\"older seminorm defined by the metric $\omega_{(\beta_{i})}$, which is of course equivalent to that defined by the metric $\omega_{i}$.
This Lemma is straightforward to prove  using the H\"ormander technique to construct a suitable section of $L^{k}\rightarrow X_{i}$ for an appropriate power $k$. The power $\rho^{2-\alpha}$ comes from the scaling behaviour: that is when the ball $B^{\beta_{i}}(q,\rho)$ is scaled to unit size and the metric is scaled the estimate for the rescaled potential $\tilde{\psi}$ becomes $[\tilde{ \psi}]_{\alpha}\leq C$. On the other hand the result is local and can probably also be proved by using the \lq\lq weight function'' version of the H\"ormander technique (with any $R>2$), or by other methods from complex analysis.

To simplify notation we drop the index $i$. We want to fix a suitable value of $\zeta$. First, consider the function $F(M)= \det M- {\rm Tr} M$ on the space of $n\times n$ Hermitian matrices. Since the derivative at the identity vanishes we can for any $\eta>0$ find a $\zeta(\eta)$ such that if $M_{1}, M_{2}$ satisfy $(1-\zeta)\leq M_{i}\leq (1+\zeta)$ then 
\begin{equation} \vert F(M_{1})-F(M_{2})\vert\leq \eta \vert M_{1}- M_{2}\vert. \end{equation}
Next we recall the Schauder estimate of \cite{kn:Dona11}. This asserts that  there is some $K=K(\alpha,\beta)$ such that for all $\psi\in {\cal C}^{2,\alpha,\beta}(B^{\beta}(0,2))$ we have
\begin{equation} [\dbd \psi]_{\alpha, B^{\beta}(0,1)}\leq K \left( [\Delta \psi]_{\alpha, B^{\beta}(0,2)} + [\psi]_{\alpha, B^{\beta}(0,2)}\right). \end{equation}
Here $\Delta$ is the Laplacian of $\omega_{(\beta)}$. More generally, we can choose $K=K(\alpha,\beta)$ such that for {\it any point} $q\in \bC^{n}$, radius $\rho>0$  and for $\psi\in {\cal C}^{2,\alpha,\beta}(B^{(\beta)}(q,2\rho))$ we have
\begin{equation} [\dbd \psi]_{\alpha, B^{\beta}(q,\rho)}\leq K \left( [\Delta \psi]_{\alpha, B^{\beta}(q,2\rho)} + \rho^{-2} [\psi]_{\alpha, B^{(\beta)}(q,2\rho)}\right). \end{equation}
This follows by a straightforward argument, using (38) (after scaling and translation) at points near the singular set and the standard Schauder estimate away from the singular set. Furthermore, for any $\beta$ and $\alpha<\beta^{-1}-1$ we can suppose that the  inequality, with a fixed $K$, holds for $\alpha',\beta'$ sufficiently close to $\alpha,\beta$. Now fix $\zeta$ so that $\zeta<\zeta_{0}$ and
\begin{equation}   (R-2)^{-\alpha} K \eta(\zeta) \leq R^{-\alpha}/4, \end{equation}
where $\zeta_{0}, R$ are as in Lemma 3. 
For distinct points $x,y$ in the Euclidean unit ball $B$, not in the singular set, define
\begin{equation}  Q(x,y)= d_{\beta}(x,y)^{-\alpha} \vert \omega(x)-\omega(y) \vert \left( \min(d_{\beta}(x,\partial B), d_{\beta}(y,\partial B)\right)^{\alpha}. \end{equation}
Here $d_{\beta}$ denotes the distance in the metric $\omega_{(\beta)}$. The difference $\omega(x)-\omega(y)$ is interpreted in the sense of \cite{kn:Dona11}, in that we take the matrix entries with respect to a fixed orthonormal frame for the (1,1) forms (in which $\omega_{(\beta)}$ is constant). Thus, by definition,
$$  [\omega]_{\alpha, B} = \sup_{x,y}\left( d_{\beta}(x,y)^{-\alpha} \vert \omega(x)-\omega(y)\vert \right), $$
where the supremum is taken over distinct points $x,y$ in $B$.
By hypothesis this supremum is finite, so we have a finite number 
$  M=\sup_{x,y} Q(x,y)$. We seek an {\it a priori} bound on $M$.  It is not immediate that this supremum is attained, but we can certainly choose $x,y$ so that $Q(x,y)\geq M/2$. Set $d=  \min(d_{\beta}(x,\partial B), d_{\beta}(y,\partial B))$.  If $d_{\beta}(x,y)\geq d/R$ 
then $M/2\leq R^{\alpha} \vert \omega(x)-\omega(y)\vert$ and we are done since we have an $L^{\infty}$ bound on $\omega$. So we can suppose that $d_{\beta}(x,y)<d/R$. Let $x$ be the point closest to the boundary and take it to be the centre $q$ in Lemma 3, with radius $\rho=d/R$. Thus we can apply Lemma 3 to write $\omega=\omega_{(\beta)} + \dbd \psi $ over $B^{\beta}(x,2\rho)$, with
$[\psi]_{\alpha} \leq C \rho^{2-\alpha}$. Thus by construction  \begin{equation}[\dbd \psi]_{\alpha, B^{\beta}(x,\rho)}\geq (M/2)d^{-\alpha}= (M/2) (R\rho)^{-\alpha},  \end{equation} since $y$ lies in $B^{\beta}(x,\rho)$. On the other hand if $z,w$ are distinct points in $B^{\beta}(x, 2\rho)$ the inequality $Q(z,w)\leq M$ implies that
$$d_{\beta}(z,w) \vert \omega(z)-\omega(w) \vert \leq ((R-2)\rho)^{-\alpha} M,$$ since the distance from $z,w$ to the boundary is at least $(R-2) \rho$. So
 \begin{equation}[\dbd \psi]_{\alpha, B^{\beta}(x,2\rho)}\leq ((R-2)\rho)^{-\alpha} M. \end{equation}
Write $\det$ for the ratio of the volume forms $\omega^{n}, \omega_{(\beta)}^{n}$.
Then (37) implies that
$$   [\det -\Delta \psi]_{\alpha, B^{\beta}(x,2\rho)}\leq \eta  [\dbd \psi ]_{\alpha, B^{\beta}(x,2\rho)}, $$
So $$ [\Delta \psi]_{\alpha, B^{\beta}(x,2\rho)}\leq \eta [\dbd \psi ]_{\alpha, B^{\beta}(x,2\rho)} + [\det]_{\alpha, B^{\beta}(x,2\rho)}\leq ((R-2) \rho)^{-\alpha} M \eta+ [\det]_{\alpha, B^{\beta}(x,2\rho)}.$$
Thus by (39) and (42)
$$ R^{-\alpha} \rho^{-\alpha} M/2\leq [\dbd \psi]_{\alpha, B^{\beta}(x,\rho)}\leq ((R-2) \rho)^{-\alpha} K M \eta+ K  [\det]_{\alpha, B^{\beta}(x,2\rho)} + K \rho^{-2} [\psi]_{\alpha, B^{\beta}(x,2\rho)}, $$
and by our choice $ (R-2)^{-\alpha} K\eta\leq R^{-\alpha}/4$ we can rearrange to get
$$ R^{-\alpha} \rho^{-\alpha} M/4\leq  K  [\det]_{\alpha, B^{\beta}(x,2\rho)} + K \rho^{-2} [\psi]_{\alpha, B^{\beta}(x,2\rho)}. $$
Now, putting in the bound from Lemma 3 and multiplying by $\rho^{\alpha}$, we have
$$  R^{-\alpha} M/4 \leq K \rho^{\alpha} [\det]_{\alpha} + C K, $$
which gives our bound on $M$. 

To establish Theorem 2 we fix $\alpha<\beta_{\infty}^{-1}-1$ and apply the above discussion to $\omega_{i}$ with $\beta=\beta_{i}$, so we can take a fixed constant $K$ in the Schauder estimate. We see that the limiting metric, as $i\rightarrow \infty$, lies in ${\cal C}^{,\alpha, \beta_{\infty}}$.
 
{\bf Remark}
To prove Theorem 2 we have only had to deal with the points in the divisor $\Delta$, since the other points are covered by arguing as in \cite{kn:DS}, using a version of the Evans-Krylov theory. On the other hand we can also handle these latter points by just the same argument used above, not invoking the Evans-Krylov theory (from PDE) but applying instead the Cheeger-Colding Theory (from Riemannian geometry).

\end{document}